\documentclass[preprint,12pt]{elsarticle}

\usepackage{graphicx} 
\usepackage{enumitem}
\usepackage[export]{adjustbox}
\usepackage{amssymb}
\usepackage{amsthm}
\usepackage{amsmath}
\usepackage{amsfonts}
\usepackage{slashbox}
\usepackage{algorithm} 
\usepackage{algorithmic}  
\usepackage[algo2e]{algorithm2e} 
\usepackage{hyperref}
\usepackage{lipsum}
\usepackage{cleveref}
\usepackage[caption=false]{subfig}
\usepackage{units}
\usepackage{lineno}

\def\r{\mathbb{R}}
\def\n{\mathbb{N}}

\def\z{\mathbb{Z}}
\def\dd{{\rm div}}
\def\trian{\mathfrak{T}} 

\def\unn{\textbf{u}}
\def\vnn{\textbf{v}}
\def\xnn{\textbf{x}}
\def\ynn{\textbf{y}}
\def\dert{\partial_t}
\def\omep{\Omega^{\varepsilon}}
\def\kep{\mathbf{K}^{\varepsilon}}
\def\knn{\mathbf{K}}
\def\kef{\mathbf{K}^{\star}}
\def\dt{\Delta t}
\def\bigo{\mathcal{O}}

\begin{document}

\begin{frontmatter}



	\title{Numerical homogenization of non-linear diffusion problems on adaptive meshes}

	\author[l1]{Manuela Bastidas\corref{cor1}}
	\cortext[cor1]{Corresponding author.}
	\ead{manuela.bastidas@uhasselt.be}

	\author[l2]{Carina Bringedal}
	\author[l1]{Iuliu Sorin Pop}
	\author[l3]{Florin Adrian Radu}

	\address[l1]{Faculty of Sciences, Hasselt University. Diepenbeek, Belgium.}
	\address[l2]{Institute for Modelling Hydraulic and Environmental Systems, University of Stuttgart. Stuttgart, Germany.}
	\address[l3]{Department of Mathematics, University of Bergen. Bergen, Norway.}

	\begin{abstract}
		We propose an efficient numerical strategy for solving non-linear diffusion problems defined in a porous medium with highly oscillatory characteristics. This scheme is based on the classical homogenization theory and uses a locally mass-conservative formulation at different scales. In addition, we discuss some properties of the proposed non-linear solvers and use an error indicator to perform a local mesh refinement. The main idea is to compute the effective parameters in such a way that the computational complexity is reduced but preserving the accuracy. We illustrate the behaviour of the homogenization scheme and of the non-linear solvers by performing two numerical tests. We consider both a quasi-periodic example and a problem involving strong heterogeneities in a non-periodic medium.
	\end{abstract}


\end{frontmatter}


\section{Introduction}
\label{sec:intro}

Non-linear parabolic problems are encountered as mathematical models for several real life applications. Examples in this sense are partially saturated flow in porous media, non-steady filtration, and reaction-diffusion systems. Realistic applications often involve heterogeneous media, which translate into highly oscillatory coefficients and non-linearities.

Letting $\Omega^{\varepsilon}$ be a bounded, possibly perforated domain in $\r^d$ ($d=2,3$) with Lipschitz boundary $\partial \Omega^{\varepsilon}$ and $\mathrm{T}>0$ be a maximal time. We consider the non-linear diffusion equation
\begin{equation}\label{eq:OriginalProblem0}
	\dert b^{\varepsilon}(\xnn,p^{\varepsilon}(\xnn,t)) - \dd\left( \kep(\xnn)\,\nabla p^{\varepsilon}(\xnn,t) \right)  = f^{\varepsilon}(\xnn,t), \quad \text{in } \Omega^{\varepsilon}\times \left( 0,\mathrm{T} \right],
\end{equation}
with suitable initial and boundary conditions. In this setting, $\varepsilon$ is a positive small parameter and denotes the scale separation between the micro-scale (e.g. the scale of pores in a porous medium) and the macro-scale (e.g. the Darcy scale, the scale of simulation in case of heterogeneous media). With the superscript $0< \varepsilon\ll 1$ we indicate that the quantities involve highly oscillatory features and the medium is considered highly heterogeneous.
Inspired by unsaturated fluid flow in a porous medium \eqref{eq:OriginalProblem0} can, for example, represent the non-dimensional Richards equation after applying the Kirchhoff transformation, without taking into account gravity effects (see \cite{bear2012introduction}). In this case, the primary unknown $p^{\varepsilon}(\xnn,t)$ is the transformation of the fluid pressure. For simplicity $p^{\varepsilon}(\xnn,t)$ will be called \textit{pressure} in what follows. The given data include the source $f^{\varepsilon}$,  the absolute permeability matrix $\kep$ and the volumetric fluid saturation $b^\varepsilon$, which is a given function of $p^{\varepsilon}$.

The key issue in developing numerical methods capturing the interaction between scales is to avoid the high computational cost. The use of classical schemes over fine-scale meshes has often unreachable requirements. To capture the oscillations in the medium the required mesh size is smaller than $\varepsilon$. In this sense, standard numerical methods will either fail or become inefficient.

There are numerous numerical simulation techniques for processes that involve two or more scales in space and time. During the last years, approaches like the multi-scale finite-volume (MSFV), the algebraic dynamic multilevel (ADM), the heterogeneous multi-scale (HMM) and the multi-scale finite element (MsFEM) methods are becoming more and more relevant.
Concretely, the MSFV and ADM methods proposed in \cite{hajibeygi2008iterative,cusini_algebraic_2016} aim to solve problems involving different scales by incorporating the fine-scale variation into the coarse-scale operators. The multi-scale finite volume method (MSFV) in \cite{cusini_algebraic_2016} includes a dynamic local grid refinement method to provide accurate and efficient simulations employing fine grids only where needed. We highlight that the MSFV and ADM use a section of the fine-scale feature to construct the macro-scale solution without estimations of the macro-scale parameters.
A preliminary study in parallel of the two approaches based on ADM and on numerical homogenization is contained in \cite{HadiManuela}.

On the other hand, the HMM (see \cite{weinan2003heterogeneous,abdulle2012heterogeneous}) relies on coupled macro and micro-scale solvers using homogenization (see \cite{hornung2012homogenization}). This method takes advantage of the scale separation and is based on the numerical approximation of the macro-scale data. In \cite{abdulle2012heterogeneous,Abdulle159424,abdulle2011adaptive} ideas on how to manage different scales in an efficient computational way are developed, using the standard finite element method (FEM). Further, the numerical computations using finite difference and discontinuous Galerkin method also demonstrate the potential of this framework in \cite{weinan2003heterogeneous,chen2005heterogeneous}.

Improved multi-scale methods to simulate non-linear single-phase and multi-phase flow has been proposed in \cite{amanbek2017adaptive, arbogast2007multiscale,wheeler2002enhanced,moyner2016multiscale}.  An Enhanced Velocity Mixed element method is proposed in \cite{wheeler2002enhanced} to deal with non-matching, multi-block grids and couple micro and macro-scale domains. In the same line of research, \cite{arbogast2007multiscale} give a computational strategy for the multi-scale dynamics over non-matching grids using mesh refinement and enriched multi-scale basis functions. In \cite{amanbek2017adaptive}, the homogenization theory is combined with domain decomposition to obtain locally effective parameters and solve macro-scale problems. Further, the  multi-scale finite element (MsFEM) method presented in \cite{hellman2016multiscale,henning2015adaptive,henning2014adaptive}
constructs a low dimensional multi-scale mixed finite element space. This change of the discrete spaces allows the formal derivation of a-posteriori estimates to control the micro-scale error and its influence on the macro-scale.

In this paper, we develop a locally mass-conservative scheme that computes the homogenized permeability field of \eqref{eq:OriginalProblem0} over coarse meshes. In contrast with the papers mentioned before, we use an error indicator on the macro-scale solvers to localize the error and subsequently refine or coarsen the mesh accordingly. We propose a combination of techniques supported in the theoretical framework of the homogenization (see \cite{hornung2012homogenization}) for non-linear parabolic equations. Our strategy relies on the solution of  micro-scale cell problems to calculate averaged parameters that are used in a macroscopic solver.
The computation of the effective parameters can be parallelized and it is cheap to perform. Moreover, the error induced by the calculation of the effective parameters can be dismissed when one applies a sufficiently accurate micro-scale solver. It is important to remark that, although periodicity is assumed in the classical homogenization theory, the upscaling technique developed here can be applied to problems in non-periodic media and we show, in the numerical examples, that the effective parameters still represents the macro-scale behaviour.

We apply the backward Euler (BE) method for the time discretization and the mixed finite element method (MFEM) for the spatial discretization. In order to solve the fully discrete formulation of \eqref{eq:OriginalProblem0}, non-linear solvers are required. We discuss  the applicability of classical iterative solvers like Newton or Picard (see \cite{bergamaschi1999mixed,celia1990general}) and we detail the formulation of a robust fixed point method called L-scheme proposed in \cite{list2016study}. This linearization procedure has the advantage of being unconditionally convergent. More exactly, the convergence of the L-scheme is neither affected by the initial guess nor by the mesh size. Nevertheless, the convergence rate of the L-scheme is only linear and therefore slower compared to the Newton scheme (see \cite{radu2017robust}).
We mention the paper \cite{list2016study} for an approach combining the L and the Newton schemes in an optimized way. There, the L-scheme is applied to provide a suitable initial point for the Newton scheme. We use this strategy to improve the convergence of the scheme up to the quadratic convergence. We also refer to \cite{mitra2018modified} for a modified L-scheme featuring improved convergence (compared to the L-scheme) and scalability properties (compared to Newton and Picard).

For time-dependent problems the idea of adaptive meshes is very useful to localize the changes in the solution between different time steps. On the other hand, reaching finer meshes becomes computationally expensive because it requires extra calculations of the macro-scale parameters. The finer the mesh for the upscaled model, the higher the computational effort as the effective parameters need to be computed in more points, thus more cell problems need to be solved. For this reason, we present an error indicator that specifies when the numerical solution and the effective parameters should be re-computed. With this strategy we aim to control the convergence rate of the numerical scheme and to avoid unnecessary computations of the local problems.

The paper is organized as follows. In \Cref{sec:model} the details of the model, the geometry and the discrete formulation are given and the necessary assumptions are stated. \Cref{sec:Periodic} gives a brief summary of the standard procedure of the homogenization for a parabolic case in a periodic porous media.  In \Cref{sec:noperiodic} the mesh refinement and the coarsening strategy is stated and in \Cref{sec:Lsche} the linearization scheme is described. We discuss the numerical tests in \Cref{sec:Num}, where the quasi-periodic and non-periodic cases are considered.

\section{The model formulation and the spatial discretization}
\label{sec:model}

To construct a robust and locally conservative scheme we consider the mixed formulation of \eqref{eq:OriginalProblem0}. Letting  $\unn^{\varepsilon}(\xnn,t)$ be the Darcy velocity, the unknowns $(p^{\varepsilon},\unn^{\varepsilon})$ satisfy
\begin{equation}\label{eq:mixedOrig}
	\begin{aligned}
		\dert b^{\varepsilon}(\xnn,p^{\varepsilon}(\xnn,t)) + \dd\left( \unn^{\varepsilon}(\xnn,t) \right) & = f^{\varepsilon}(\xnn,t),                       &  & \text{in } \omep_{\mathrm{T}},          \\
		\unn^{\varepsilon}(\xnn,t)                                                                         & = - \kep(\xnn)\, \nabla p^{\varepsilon}(\xnn,t), &  & \text{in } \omep_{\mathrm{T}},          \\
		p^{\varepsilon}(\xnn,t)                                                                            & = 0 ,                                            &  & \text{on } \partial \omep_{\mathrm{T}}, \\
		p^{\varepsilon}(\xnn,0)                                                                            & = p_I,                                           &  & \text{in }  \omep,
	\end{aligned}
\end{equation}

Here $\omep_{\mathrm{T}} := \omep \times (0,\mathrm{T}]$ and $\partial \Omega^{\varepsilon}_\mathrm{T} := \partial \omep \times (0,\mathrm{T}]$. As mentioned before, by using the superscript $\varepsilon>0$ we emphasize on the fact that rapidly oscillating characteristics are involved. For example, the domain either involves characteristics changing within $\varepsilon$-sized regions, or it may include perforations.

We refer to \cite{alt1983quasilinear} for the existence and uniqueness of a weak solution of \eqref{eq:OriginalProblem0} under the following assumptions:
\crefname{enumi}{part}{assumptions}
\begin{enumerate} [label={(A\arabic*)}]
	\item \label{A1}  The function $b^{\varepsilon}(\xnn,\cdot)$ is non-decreasing, $b^{\varepsilon}(\cdot,0) =0$ and
	      H\"older continuous.  There exists $\alpha \in (0,1]$ and $L_b > 0$ such that
	      \[| b^{\varepsilon}(\xnn,p_1) -  b^{\varepsilon}(\xnn,p_2) | \leq L_b |p_1-p_2|^{\alpha},\]
	      for all  $\xnn \in \omep$ and $p_1,\,p_2 \in \r$.
	\item \label{A3} The permeability function $\kep:\Omega^{\varepsilon}  \to \r^{d\times d}$ is symmetric for all $\xnn \in \Omega^{\varepsilon}$ and continuous. There exist $\beta,\lambda>0$ such that
	      \begin{equation*}
		      \beta \|\psi\|^2 \leq \psi^\mathsf{T} \, \kep(\xnn) \, \psi  \leq \lambda \|\psi\|^2  \quad \text{ for all } \psi \in \r^d \text{ and } \xnn \in \omep.
	      \end{equation*}
	\item \label{A4}The initial data $p_I$ and the source term $f^{\varepsilon}$ are  essentially bounded uniformly w.r.t $\varepsilon$.
\end{enumerate}
In \cite{radu2008error} the equivalence between the mixed and conformal weak formulations is proved in both continuous and semi-discrete cases.

\subsection{The non-linear fully discrete problem}
\label{sec:discOrig}

To define the discrete problem we let $\trian_{h^{\varepsilon}} $ be a triangular partition of the domain $\omep$ with elements $\mathcal{T}$ of diameter $h^{\varepsilon}_{\mathcal{T}}$ and $h^{\varepsilon} := \max\limits_{\mathcal{T}\in\trian_{h^{\varepsilon}}}{h^{\varepsilon}_{\mathcal{T}}}$ such that $h^{\varepsilon}\ll\varepsilon$.

Further, $0= t_0 \le t_1 \le t_1 \le \dots \le t_N = \mathrm{T}$, $N \in \n$ is a partition of the time interval $[0,\mathrm{T}]$ with constant step size $\dt = t_{i+1}-t_i$, $i\geq 0$. For the discretization of the flux $\unn^{\varepsilon}$ we consider the lowest-order Raviart-Thomas space $V_{h^{\varepsilon}}:=\mathcal{R}T_0(\trian_{h^{\varepsilon}})$ and  for the pressure $p^\epsilon$ we use the discrete subspace of piecewise constant functions $W_{h^{\varepsilon}}$ (see \cite{brezzi2012mixed})
\begin{equation*}
	\begin{aligned}
		W_{h^{\varepsilon}} & := \left\lbrace q \in L^2(\omep) \, | \, q \text{ is constant on each element } \mathcal{T} \in \trian_{h^{\varepsilon}} \right\rbrace,                                                         \\
		V_{h^{\varepsilon}} & := \left\lbrace \vnn \in H(\dd,\omep) \, | \, \vnn|_{\mathcal{T}} = \mathbf{a}+b\xnn \text{ for all } \mathcal{T} \in \trian_{h^{\varepsilon}},\, \mathbf{a} \in \r^d,\, b\in \r \right\rbrace,
	\end{aligned}
\end{equation*}
with $L^2(\omep)$ being the space of the square-integrable functions with the usual norm and $H(\dd,\omep) :=  \left\lbrace \vnn \in [L^2(\omep)]^d \,|\, \dd(\vnn) \in L^2(\omep) \right\rbrace$. We let $\left\langle \cdot, \cdot \right\rangle$ represent the inner product on $L^2(\omep)$.

\paragraph{\textbf{Problem} $\mathbf{PM^\varepsilon_n}$}\label{problemPMn}
Let $n \geq 1$. Given $(\left( p^{\varepsilon}\right)_{h^{\varepsilon}}^{n-1},\left( \unn^{\varepsilon}\right)_{h^{\varepsilon}}^{n-1}) \in W_{h^{\varepsilon}}\times V_{h^{\varepsilon}}$, find $\left( p^{\varepsilon}\right)_{h^{\varepsilon}}^n \in W_{h^{\varepsilon}}$ and $\left( \unn^{\varepsilon}\right)_{h^{\varepsilon}}^n \in V_{h^{\varepsilon}}$ such that for any $q \in W_{h^{\varepsilon}}$ and $\vnn \in V_{h^{\varepsilon}}$ there holds
\begin{equation*}\label{eq:discretemixedori} 
	\begin{aligned}
		\left\langle b^{\varepsilon}\left( \cdot,\left( p^{\varepsilon}\right)_{h^{\varepsilon}}^n \right)  -  b^{\varepsilon}\left( \cdot,\left( p^{\varepsilon}\right)_{h^{\varepsilon}}^{n-1}\right) , q \right\rangle  + \dt \, \left\langle \dd\left( \left( \unn^{\varepsilon}\right)_{h^{\varepsilon}}^n \right), q\right\rangle & =  \quad \dt\, \left\langle f^{\varepsilon}, q \right\rangle, \\
		\left\langle  \left[ \kep \right]^{-1} \,  \left( \unn^{\varepsilon}\right)_{h^{\varepsilon}}^n, \vnn \right\rangle  - \left\langle   \left( p^{\varepsilon}\right)_{h^{\varepsilon}}^n,  \dd\left(  \vnn \right)  \right\rangle                                                                                                & = \quad 0.
	\end{aligned}
\end{equation*}
We denote by $\left( p^{\varepsilon}\right)_{h^{\varepsilon}}^0$ the $L^2$-projection of the initial condition $p_I$ over the mesh $\trian_{h^{\varepsilon}}$. For simplicity, we omit writing the $\xnn$ argument in $ b^{\varepsilon}(\xnn,p^\varepsilon)$, which becomes now  $b^{\varepsilon}(p^{\varepsilon})$.

For details about the existence and uniqueness of the solution to problem $\mathbf{PM^\varepsilon_n}$ we refer to \cite{radu2008error}. Note that the problem $\mathbf{PM^\varepsilon_n}$ is non-linear, therefore a non-linear solver is needed. This is detailed in \Cref{sec:Lsche}.

\section{The two-scale approach}
\label{sec:Periodic}

We start the presentation for the case of a periodic medium. Building on this, we extend these ideas for non-periodic situations. The concept of coupling the scales trough the calculation of effective parameters is used, among others, in \cite{abdulle2012heterogeneous,abdulle2011adaptive,amanbek2017adaptive}. Here we follow the ideas therein and enhance the strategy with adaptive mesh refinement and robust non-linear solvers.

\begin{figure}[htpb]
	\centering
	\includegraphics[width=0.7\textwidth,trim={0cm 0cm 0cm 0cm},clip]{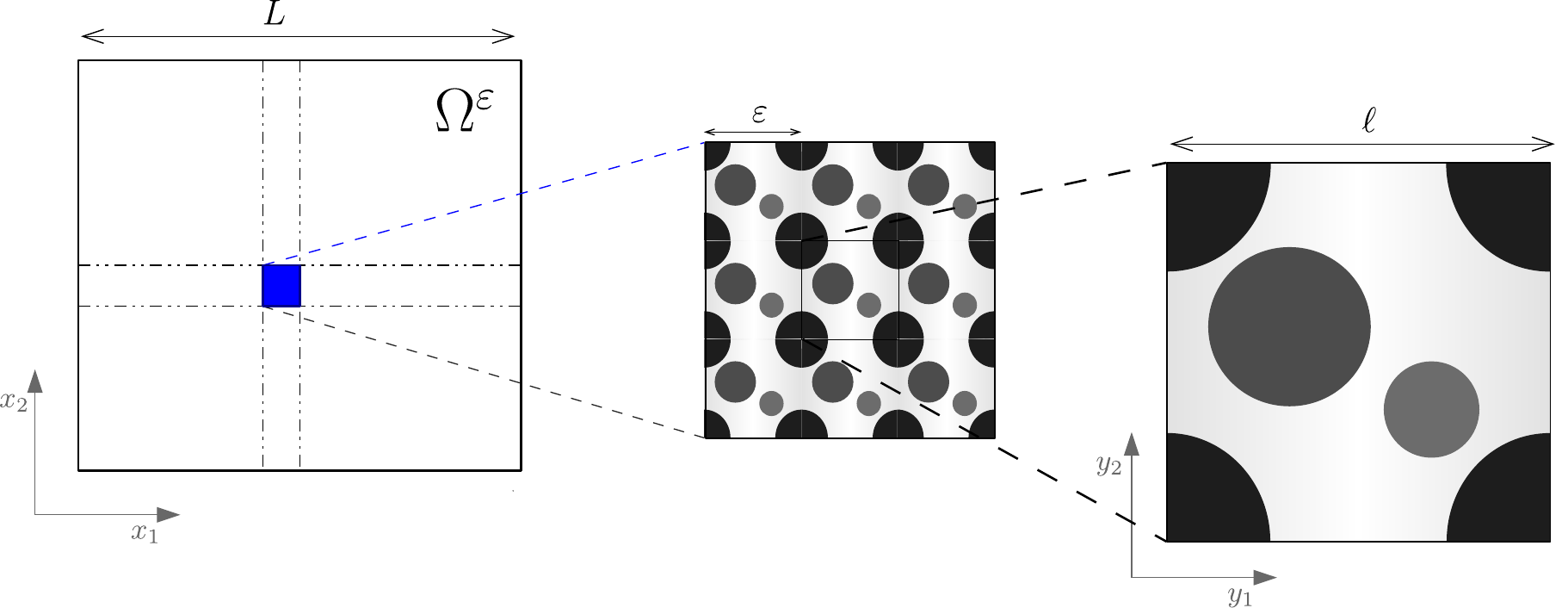}
	\vspace{-0.3cm}
	\caption{Two-scale structure. Zoom in to the pore structure in $\r^2$ where typical length sizes are indicated.}\label{fig:upscalingFigure}
\end{figure}

We assume that the domain $\omep$ can be written as the finite union of \textit{micro-scale} regions, namely $Y$, where the parameters change rapidly. In other words, the parameters and non-linearities take different values inside of $Y$ (see \Cref{fig:upscalingFigure}).
In the extreme case, the micro-scale $Y$ can be viewed as a perforated region with a pore space and a solid grain (see e.g \cite{hornung2012homogenization}). Here we give the ideas for non-perforated domains but this can be adapted straightforwardly to perforated ones.

At the micro-scale $Y$ and the macro-scale $\Omega^{\varepsilon}$ we assume characteristic lengths $\ell$ and $L$ respectively. The factor $\varepsilon:= \frac{\ell}{L}$ denotes the scale separation between the two scales. To identify the variations at the micro-scale we define a fast variable $\ynn:= \frac{\xnn}{\varepsilon}$. To each macro-scale point $\xnn \in \omep$ corresponds one micro-scale cell $Y$ that captures the fast changes in the parameters.

In the non-dimensional setting, the local cells are $Y:=[0,1]^d$ and we let $\vec{i} \in \z^d$ and $\omep = \cup \left\lbrace \varepsilon(\vec{i}+Y)\,|\,\vec{i} \in \mathcal{I}_{\varepsilon} \right\rbrace$ for some set of vector indices $\mathcal{I}_{\varepsilon}$.

To formulate the homogenized problem, we make the following assumptions:
\crefname{enumi}{part}{assumptions}
\begin{enumerate} [label={(B\arabic*)}]
	\item \label{B1} There exists a function $b:\Omega^{\varepsilon} \times \r^d \times \r \to \r$ such that $b^{\varepsilon}(\xnn,p^{\varepsilon}) := b(\xnn,\frac{\xnn}{\varepsilon},p^{\varepsilon})$ and $ b(\xnn,\cdot,p^{\varepsilon})$ is $Y$-periodic.

	\item \label{B2} There exists a function $\knn:\Omega^{\varepsilon} \times  \r^d \to \r^{d\times d}$ such that $\kep(\xnn) := \knn(\xnn,\frac{\xnn}{\varepsilon})$ where $\knn(\xnn,\ynn)$ is symmetric and continuous for all $(\xnn,\ynn) \in \Omega^{\varepsilon}\times Y$ and $\knn(\xnn,\cdot)$ is $Y$-periodic.

\end{enumerate}

\subsection{The homogenization approach}
\label{sec:Homog}

A direct numerical approximation of the problem $\mathbf{PM^\varepsilon_n}$ requires the usage of an extremely fine mesh to capture all the changes in the characteristics of the medium. We consider a homogenization-based approach and compute an effective model involving only the essential variations of the permeability matrix.
Alternative approaches, like the harmonic average of the permeability, are broadly used in \cite{renard1997calculating,lie2019introduction}. Nevertheless, those ideas are rather suited for particular cases, e.g layered media (see \cite{lie2019introduction}).

We restrict the presentation to the minimum needed for explaining the approach. We make use of the \textit{homogenization ansatz} and refer to \cite{hornung2012homogenization,cioranescu_introduction_1999} for a detailed presentation of the method.

First, we assume that $p^{\varepsilon}$ can be formally expanded as
\begin{equation}\label{eq:expansionp}
	p^\varepsilon(\xnn,t) = p(\xnn,t) + 	\varepsilon p_1(\xnn,\ynn,t) + \varepsilon^2 p_2(\xnn,\ynn,t) + ...\, ,
\end{equation}
where $\ynn=\frac{\xnn}{\varepsilon}$ stands for the fast variable, $\xnn$ is the slow variable and each function $p_i:\omep\times Y\times (0,T] \to \r$ is $Y$-periodic w.r.t $\ynn$. The function  $p(\xnn,t)$ does not depend on $\ynn$ and is in fact the \textit{macro-scale approximation} of the pressure $p^{\varepsilon}(\xnn,t)$.

Additionally, the two-scale gradient and divergence operators become
\begin{equation}\label{twoscaleop}
	\nabla = \nabla_x +\frac{1}{\varepsilon}\nabla_y \quad \text{and} \quad \dd = \dd_x+\frac{1}{\varepsilon}\dd_y.
\end{equation}

Using \eqref{eq:expansionp} and \eqref{twoscaleop} in \eqref{eq:mixedOrig} and applying the Taylor expansion of $b(\cdot,\cdot,p)$ we obtain
\begin{equation*}
	\dert b -  \left( \dd_x+\frac{1}{\varepsilon}\dd_y \right) \left( \knn\left( \nabla_x +\frac{1}{\varepsilon}\nabla_y\right) \left( p +\varepsilon p_1 + \varepsilon^2p_2  \right) \right)  + \bigo\left({\varepsilon}\right)  = f.
\end{equation*}

To determine $p_1$ as a function of $p$, for the terms of order $\bigo(\varepsilon^{-1})$ we can write
$p_1(\xnn,\ynn,t)= \hat{p_1}(\xnn,t)+\sum_{j=1}^{d} \frac{\partial p(\xnn,t)}{\partial x_j} \, \omega^j(\xnn,\ynn)$
where the function $\hat{p_1}$ is an arbitrary function of $\xnn$, and $\omega^j$ are the $Y$-periodic solutions of the following mixed \textit{micro-cell} problems
\begin{equation}\label{eq:cellprob}
	\begin{aligned}
		\dd_y  \boldsymbol{\xi}^j & =   \dd_y \left(  \knn(\xnn,\cdot)\, \mathbf{e}_j  \right), & \text{ in } Y, \\
		\boldsymbol{\xi}^j        & =  - \knn(\xnn,\cdot)\,\nabla_y \omega^j,                   & \text{ in } Y.
	\end{aligned}
\end{equation}
Here $\{\mathbf{e}_j\}_{j=1}^{d}$ is the canonical basis of $\r^d$. To guarantee the uniqueness of the solution we assume that $\omega^j$ has the average $0$ over the micro cells, that is, $\int_Y \omega^j(\xnn,\ynn)d\ynn = 0$ for all $\xnn\in \omep$.

To simplify the notation, we use $\Omega$ instead of $\Omega^{\varepsilon}$ for the macro-scale domain and  $\partial\Omega$ for its outer boundary. Following from the homogenization, $\Omega$ does not contain any oscillatory behaviour. Recalling the periodic boundary conditions and averaging over $Y$ one obtains the \textit{homogenized} mixed formulation.

Letting $\unn(\xnn,t)$ denote the upscaled Darcy velocity, the upscaled unknowns $(p,\unn)$ satisfy
\begin{equation}\label{eq:mixedUps}
	\begin{aligned}
		\dert b^{\star}(\xnn,p(\xnn,t)) + \dd\left( \unn(\xnn,t) \right) & = f^{\star}(\xnn,t),               &  & \text{in } \Omega_{\mathrm{T}}           \\
		\unn(\xnn,t)                                                     & = - \kef(\xnn)\, \nabla p(\xnn,t), &  & \text{in } \Omega_{\mathrm{T}},          \\
		p(\xnn,t)                                                        & = 0 ,                              &  & \text{on } \partial \Omega_{\mathrm{T}}, \\
		p(\xnn,0)                                                        & = p_I,                             &  & \text{in }  \Omega.
	\end{aligned}
\end{equation}
Here $\Omega_{\mathrm{T}}:= \Omega\times (0,{\mathrm{T}}]$ and $\partial \Omega_{\mathrm{T}}:= \partial \Omega \times (0,\mathrm{T}]$. The effective permeability $\kef:\Omega\to \r^{d\times d}$ has the elements
\begin{equation}\label{eq:EffectiveTensor}
	\kef_{i,j}(\xnn) = \int_Y \left(  \knn(\xnn,\ynn) \left( \mathbf{e}_j + \nabla_y \omega^j(\xnn,\ynn)\right) \right) \cdot \mathbf{e}_i \, d\ynn, \quad (i,j=1,\dots, d).
\end{equation}

The upscaled saturation and source terms are
\begin{equation*}\label{eq:Effectivebf}
	b^{\star}(\xnn,p)  :=  \int_Y b(\xnn,\ynn,p) \, d\ynn \quad \text{ and } \quad
	f^{\star}(\xnn,t)  :=  \int_Y f(\xnn,\ynn,t) \, d\ynn.
\end{equation*}

The difference between the solution of \eqref{eq:mixedOrig} and the solution of \eqref{eq:mixedUps} is subtle. In the original problem, the main characteristics are present at all scales in a strongly coupled manner. The \textit{homogenized model} instead involves only essential variations at the macro-scale. However, to determine the value of the permeability tensor at a macro point $\xnn \in \Omega$, one has to solve $d$ micro-cell problems \eqref{eq:cellprob} associated with that macro point.  Note that these problems reflect the rapidly oscillating characteristics and are decoupled from the macro-scale variations. From a computational point of view, the importance of this decoupling becomes obvious. Instead of solving the full problem on a very fine mesh, one solves a collection of simpler problems. In general, analytic solutions are not available to compute the homogenized parameters. Then $\kef$, $b^{\star}$ and $f^{\star}$ must usually be computed numerically and can therefore only be obtained at discrete points of the domain $\Omega$.

If the original permeability $\kep$ satisfies \ref{A3} and \ref{B2} then the effective tensor in \eqref{eq:EffectiveTensor} is also symmetric and positive definite. Nevertheless, the numerical approximation to the effective tensor can contain non-zero non-diagonal components or different diagonal components. To quantify the anisotropic deviation of $\kef$ we compute the following quantities
\begin{equation*}
	\tau_1 = \left( \frac{\int_\Omega ||\kef_{D}(\xnn)-\kef(\xnn)||_2^2 d\xnn}{\int_\Omega ||\kef_{D}(\xnn)||_2^2 d\xnn}\right)^{\frac{1}{2}} \text{ and } \tau_2 = \left( \frac{\int_\Omega  |\kef_{1,1}(\xnn)-\kef_{2,2}(\xnn)|^2 d\xnn}{\int_\Omega  \frac{\kef_{1,1}(\xnn)^2}{2}+\frac{\kef_{2,2}(\xnn)^2}{2}d\xnn} \right)^{\frac{1}{2}}
\end{equation*}
where $\kef_{D}$ is the diagonal matrix that contains the diagonal elements of $\kef$.
The non-linear discrete problem associated with the homogenized formulation \eqref{eq:mixedUps} is defined in the following sections.

\subsection{The non-linear fully discrete homogenized problem}
\label{secsub:nonlinear}

Let $\trian_H $ be a coarse, triangular partition of the domain $\Omega$ with coarse elements $\mathcal{T}$ of diameter $H_\mathcal{T}$ and $H := \max\limits_{\mathcal{T}\in\trian_H}H_\mathcal{T}$. For the discretization of the flux $\unn$ we consider the lowest-order Raviart-Thomas space $V_H:=\mathcal{R}T_0(\trian_H)$ and  for the pressure $p$ we use the discrete subspace of piecewise constant functions $W_H$ (see \cite{brezzi2012mixed}).
\paragraph{\textbf{Problem} $\mathbf{PH_n}$}\label{problemPHn}
For a given $p_H^{n-1} \in W_H$ and $n\geq 1$, find $p_H^n \in W_H$ and $\unn_H^n \in V_H$ such that for any $q_H \in W_H$ and $\vnn_H \in V_H$ there holds
\begin{equation*}\label{eq:discretemixed} 
	\begin{aligned}
		\left\langle b^{\star}\left( \cdot,p_H^n \right)  -  b^{\star}\left(\cdot, p_H^{n-1}\right) , q_H \right\rangle  + \dt \, \left\langle \dd\left( \unn_H^n \right), q_H \right\rangle & =  \quad \dt\, \left\langle f^{\star}, q_H \right\rangle, \\
		\left\langle  \left[ \kef \right]^{-1} \,  \unn_H^n, \vnn_H \right\rangle  - \left\langle   p_H^n,  \dd\left(  \vnn_H \right)  \right\rangle                                         & = \quad 0.
	\end{aligned}
\end{equation*}
Again $p_H^0$ is the $L^2$-projection of the initial $p_I$ over the coarse mesh $\trian_H$. For simplicity, we omit writing the $\xnn$ argument in $ b^{\star}(\xnn,p)$, which becomes now  $b^{\star}(p)$.

\subsection{The micro-cell problems and the micro-scale discretization}
As mentioned before, the effective parameters must be computed at each integration point on the coarse triangulation $\trian_H$. The effective tensor $\kef$ depends on the solution of the micro cell problems \eqref{eq:cellprob}. To solve \eqref{eq:cellprob} we use the same MFEM scheme as for \eqref{eq:mixedUps}.

To approximate the solution of \eqref{eq:cellprob} we use a triangular decomposition $\trian_h$ of the micro-scale domain $Y$ with micro-scale mesh size $h$. For the discretization of the micro-scale unknowns we consider the lowest-order Raviart-Thomas space $V_h:=\mathcal{R}T_0(\trian_h)$ and the discrete subspace of piecewise constant functions $W_h$. At each integration point $\xnn \in \mathcal{T}$ with $ \mathcal{T} \in \trian_H$, the discrete micro-cell problem is

\paragraph{\textbf{Problem} $\mathbf{Ph_j}$}\label{problemPhj}
Find $(\omega_h^j,\boldsymbol{\xi}_h^j)\in W_h\times V_h$ satisfying
\begin{equation*}\label{eq:discMicrocell}
	\begin{aligned}
		\left\langle \dd \boldsymbol{\xi}_h^j, q_h \right\rangle                                                                                                                       & =  \quad  \left\langle \nabla\cdot \left( \mathbf{K}(\xnn,\cdot) \textbf{e}_j\right), q_h \right\rangle \\
		\left\langle  \left[ \mathbf{K}(\xnn,\cdot)\right]^{-1} \,  \boldsymbol{\xi}_h^j, \vnn_h \right\rangle  - \left\langle   \omega_h^j,  \dd\left(  \vnn_h \right)  \right\rangle & = \quad 0 ,                                                                                             \\
		\omega_h^j \text{ is } Y-\text{periodic},
	\end{aligned}
\end{equation*}
for all $q_h \in W_h$, $\vnn_h \in V_h$ and $j=1,\dots,d$.

After solving the problems $\mathbf{Ph_j}$, we use \eqref{eq:EffectiveTensor} to compute the discrete effective permeability and solve the discrete problem $\mathbf{PH_n}$. The cell problems $\mathbf{Ph_j}$ are linear problems that only need to be solved initially, or when the mesh changes. The numerical cost of solving the micro-scale problems is minor compared to the one of solving the original problem.

\subsection{Non-periodic case}
Until now the two-scale approach has been referenced by assuming periodicity of the permeability $\kep$. Nevertheless, we claim that the same strategy can be applied for non-periodic structures.
When the permeability field $\kep$ is non-periodic, the periodic boundary conditions in the problems $\mathbf{Ph_j}$ are artificially imposed. However, the problems $\mathbf{Ph_j}$ are well defined and will yield to one upscaled tensor $\kef$. In other words, when one solves the micro-cell problems the resulting effective permeability field can systematically be considered an upscaled quantity obtained from the original data.
The main issue is whether this upscaled permeability reflects the effective behaviour at the macro-scale. Hence, we combine the numerical homogenization with mesh adaptivity to capture the local variability.
In the numerical examples we show that the adaptive numerical homogenization applied to the non-periodic cases produce profitable results.

\section{The two-scale discretization}
\label{sec:noperiodic}

\begin{figure}[H]
	\centering
	\includegraphics[width=.6\textwidth]{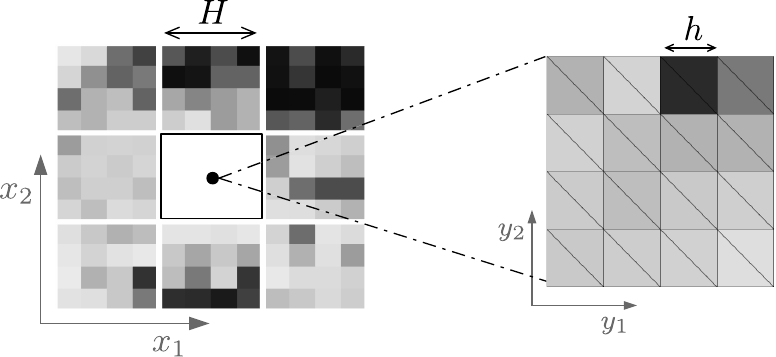}
	\vspace{-0.3cm}
	\caption{Sketch of the macro-scale partition and the correspondent micro-scale discretization in a domain $\Omega \subset \r^2$. Different intensities represent different values of the permeability.}\label{fig:Sketchparti}
\end{figure}

In practical cases, one does not necessarily have any structure in the oscillations of the data. Nevertheless the computation of macro-scale parameters remains a suitable idea. We propose to solve the micro-cell problems $\mathbf{Ph_j}$ and compute the macro-scale parameters over a coarse mesh defined beforehand. This procedure consists in two steps:
\begin{itemize}
	\item  \textit{The macro-scale partition}: Define a macro-scale division of the domain $\Omega$ with elements $Q_k$, ($k=1,2,\dots,M$), where $M$ is the total number of coarse cells.
	\item  \textit{The micro-scale domains}: Solve the micro-cell problems $\mathbf{Ph_j}$ over each coarse element $Q_k$. Note that $Q_k$ determines a micro-scale domain and there we define a micro-scale mesh size $h$.
\end{itemize}
Subsequently one can mesh the macro-scale domain and solve the homogenized problem $\mathbf{PH_n}$. In \Cref{fig:Sketchparti} we show the configuration of the macro and micro-scale partition and the procedure described previously.

\subsection{The error indicator}
\label{sec:Adapti1}

We propose a three-step strategy to adapt the macro-scale mesh to the evolution of the numerical solution of the homogenized problem. Our strategy is based on the idea of \textit{error control based on averaging technique} introduced in \cite{carstensen2004all,Carstensen99aposteriori}. The indicator of error uses a smoother approximation to the discrete solution $\unn^n_H$. We define an average operator $\mathfrak{A}_z$
\[\mathfrak{A}\unn^n_H(z) = \mathfrak{A}_z(\unn^n_H):= \frac{1}{|w_z|} \int_{w_z} \unn^n_H \, d\xnn \] where $w_z := \text{int}\left( \cup \left\lbrace \mathcal{K}\in \trian_{H_{n}}: \mathcal{K}\cap\mathcal{T}\neq \emptyset, z\in \mathcal{T} \right\rbrace \right)$ is the patch corresponding to the point $z\in\Omega$.

\subsection{The macro-scale mesh refinement}
\label{sec:Adapti2}

Our approach consists of the sequence: Solve - select the cells/triangles - refine/coarse the mesh. The mesh refining generates a sequence of  triangular meshes (one mesh per time step).

\crefname{enumi}{part}{steps}
\begin{enumerate} [label=(S\arabic*)]
	\item \label{C1} \textit{Solve}:
	      The starting point is an initial coarse mesh $\trian_{H_0}$ and the approximation of the pressure and velocity $(p_H^{0},\unn_H^{0})$ that satisfy the discrete problem $\mathbf{PH_n}$ in the first time step.

	\item \label{C2}\textit{Select the cells/triangles}:
	      Let the solution $(p_H^{n},\unn_H^{n})$ over $\trian_{H_{n}}$ be given. Calculate the error indicator $\left( \eta^n_{\mathcal{T}}\right):= \| \unn^n_H -  \mathfrak{A}\unn^n_H \|_{L^2(\mathcal{T})}$ for all $\mathcal{T} \in\trian_{H_{n}}$. The elements marked to be refined are $\mathcal{T} \in\trian_{H_{n}}$ such that (see \cite{carstensen_error_2006}) \[\eta^n_{\mathcal{T}} \geq \Theta_r\left( \max_{\mathcal{K}\in\trian_{H_{n}}} \eta^n_{\mathcal{K}} \right)  \quad \text{with } \Theta_r \in(0,1).\]
	      On the other hand, we select a set of triangles to be coarsened, i.e $\mathcal{T} \in\trian_{H_{n}}$ such that  \[\eta^n_{\mathcal{T}} \leq \Theta_c\left( \min_{\mathcal{K}\in\trian_{H_{n}}} \eta^n_{\mathcal{K}}\right)    \quad \text{with } \Theta_c\geq1.\]

	\item \label{C4} \textit{Adapt the mesh}:
	      The last step of the adaptive procedure consists of including new elements, deleting the elements to be coarsened and re-meshing. Our strategy avoids nonconforming meshes. We refine each selected cell into four new cells to compute four new effective permeabilities, and the reverse process when coarsening is necessary. Inside of the new finer cells we re-mesh with the necessary triangles.
\end{enumerate}

\begin{figure}[H]
	\centering
	\includegraphics[width=.75\textwidth]{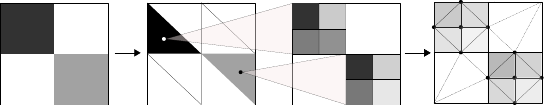}
	\caption{Outline of the mesh refinement in $\r^2$. (Left to right) Initial effective permeability. Initial triangulation and selected triangles to refine. Refinement of the permeability field. Refinement of the triangular mesh.}\label{fig:SketchRef}
\end{figure}

\begin{figure}[H]
	\centering
	\includegraphics[width=.75\textwidth]{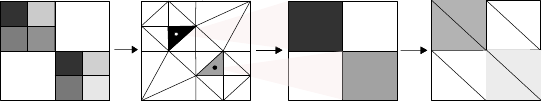}
	\caption{Outline of the mesh coarsening in $\r^2$. (Left to right) Refined effective permeability. Refined triangulation and selected triangles to coarsen. Coarsened permeability field.  Coarsened triangulation.}\label{fig:SketchRef0}
\end{figure}

The outline of the \cref{C1,C2,C4} is presented in \Cref{fig:SketchRef,fig:SketchRef0} for the 2D case. In  \Cref{fig:SketchRef} we sketch the situation when only refinement is encountered and in \Cref{fig:SketchRef0} we sketch the coarsening process.
We will only consider 2D numerical examples, but in 3D the mesh refinement can be done as described in \cite{GOLIAS1997331}.
In \Cref{fig:SketchRef,fig:SketchRef0} we highlight that at every time step it is necessary to  ensure that in the new mesh each element corresponds only to one permeability value. That restriction forces us to refine/coarsen also neighbouring elements.


With this strategy we allow to have more than one level of refinement, although the homogenization theory presented here is restricted to only two levels. Moreover, the thresholds for the refinement can be chosen depending on the problem and higher values of $\Theta_r$ and $\Theta_c$ lead to coarser meshes and less error control. We remark that the error indicator can easily be changed and the adaptive upscaling strategy does not have to be modified.

\section{The linearization scheme and the final algorithm}
\label{sec:Lsche}

A popular strategy to solve non-linear problems is the Newton method (see \cite{bergamaschi1999mixed}). The reason to use the Newton method is the quadratic convergence, but we remark that quadratic convergence only arises under certain restrictions. Specifically, the initial guess for the iterations must to be close enough to the expected solution for the scheme to converge.
For evolution equations, this means that the time step should be small enough and this usually leads to impractical values.
For that reason, we use the L-scheme which is a fixed point iteration scheme. Although it is only linearly convergent, the convergence is guaranteed regardless of the initial guess and it does not involve any derivatives (see \cite{list2016study,pop2004mixed,slodicka2002robust}).

Let $\mathfrak{L}\geq \max_{p\in\r}\left\lbrace \partial_p b^\star(\cdot,p) \right\rbrace$ be fixed and assume $p_H^{n-1}$ given. With $i\in \n$, $i\geq1$ being the iteration index, the next iteration in the L-scheme is the solution of the following linear problem.

\paragraph{\textbf{Problem} $\mathbf{PH_n^i}$}
\textit{Find $p_H^{n,(i)} \in W_H$ and $\unn_H^{n,(i)} \in V_H$ such that for any $q_H \in W_H$ and $\vnn_H \in V_H$ there holds}
\begin{equation*}\label{eq:Lscheme}
	\begin{aligned}
		\left\langle \mathfrak{L}\left( p_H^{n,(i)}- p_H^{n,(i-1)} \right) + b^{\star}\left(\cdot, p_H^{n,(i-1)}\right)  , q_H \right\rangle                                                                                                                                               \\
		+ \dt \, \left\langle \dd\left( \unn_H^{n,(i)} \right), q_H \right\rangle                                                                                     & =   \dt\, \left\langle f^{\star}, q_H \right\rangle + \left\langle  b^{\star}(\cdot,p_H^{n-1}), q_H \right\rangle, \\
		\left\langle  \unn_H^{n,(i)}, \vnn_H \right\rangle                             - \left\langle  \kef\,  p_H^{n,(i)},  \dd\left(  \vnn_H \right)  \right\rangle & =  0.
	\end{aligned}
\end{equation*}
The natural choice for the initial iteration $p_H^{n,(0)}$ is $p_H^{n-1}$. In the non-linear solver the iterations take place until one reaches a prescribed threshold for the $L^2$-norm of the residual $\partial p_H^{n,(i)} :=  p_H^{n,(i)} -   p_H^{n,(i-1)}$.

The use of an upper bound of $\partial_p b^{\star}(\cdot,p)$ affects the convergence rate. For the L-scheme the convergence rate is  $\alpha = \frac{\mathfrak{L}-m}{\mathfrak{L}+C\dt}$ for some $C>0$ and $m<\mathfrak{L}$ (see \cite{pop2004mixed}). This can lead to non-optimal convergence, for example when $\mathfrak{L}$ is very large. For this reason in \Cref{sec:Num} we choose a smaller value $\mathfrak{L}= \frac{1}{2} \max_{p\in\r}\left\lbrace \partial_p b^{\star}(\cdot,p) \right\rbrace$  which still gives convergence (see \cite{list2016study}).

Finally, we combine the non-linear solver, the mesh adaptivity and the homogenization ideas in a simple algorithm.

\begin{algorithm}[H]
	\SetAlgoLined
	\KwResult{Pressure $p_{H_N}$ and velocity $\unn_{H_N}$ over a refined mesh $\trian_{H_{N}}$}
	Choose an initial coarse-mesh $\trian_{H_{0}}$ and compute the coarse effective permeability $\kef$.

	\For{ time step $t^n$}{
	Estimate the error of the solution $\unn_{H_{n-1}}$\;

	Refine/coarsen the mesh $\trian_{H_{n-1}}$\;

	\If{new/deleted elements}{
		Re-compute the effective parameter $\kef$\;
	}\;

	\While{ $\|\partial p_H^{n,(i)}\| > tol$  }{
	Compute the solutions $p_{H_n}^{n,(i)}$ and $\unn_{H_n}^{n,(i)}$ over the new mesh $\trian_{H_{n}}$\;
	}\;	}\;
	\caption{Adaptive numerical homogenization}
\end{algorithm}

\section{Numerical results}
\label{sec:Num}

We present two numerical examples in $\r^2$ to illustrate the behaviour of the proposed adaptive homogenization procedure. We first verify our numerical homogenization approach using a manufactured periodic and quasi-periodic media and subsequently use a non-periodic test case. Note that all parameters specified in the following examples are non-dimensional and the pressures are also shifted to lie between 0 and 1.

\subsection{The periodic and quasi-periodic cases}
\label{sec:Num1}

Consider the macro-scale domain $\Omega = [0,1]\times[0,\frac{1}{2}]$ with initial condition  $p_0 = 0$ and no-flux boundary conditions. The isotropic periodic permeability field is defined by
\begin{equation*}
	\kep(\xnn) = \left(  10x_1^2x_2+\frac{1}{2+1.8\cos(2\pi\frac{x_1}{\varepsilon})\cos(2\pi\frac{x_2}{\varepsilon})}\right) \mathbb{I}_{2\times2}
\end{equation*}
A source and a sink are placed in the upper-right and the lower-left corners, having fixed pressures of $1$ and $0$, respectively. The volumetric concentration is $b^{\varepsilon}(\xnn,p) = \mathcal{R}\cdot(p^{\varepsilon})^3$. Here $\mathcal{R}$ is a non-dimensional constant that let us simulate a \textit{fast diffusion} process.
For the time discretization we take $\mathrm{T}=1$ with $\dt=0.02$.

To solve the problem $\mathbf{PM^\varepsilon_n}$ with the necessary resolution to capture the oscillations over $\Omega$ the mesh size is restricted to be $h^{\varepsilon} \ll \varepsilon$. We use $h^{\varepsilon}= 5\times10^{-3}$ to compute the fine-scale solutions $(p_{h^{\varepsilon}},\unn_{h^{\varepsilon}})$ when $\varepsilon=\frac{1}{8},\frac{1}{16}$ and $\frac{1}{32}$. The reference solutions are computed using the same MFEM, backward Euler scheme and the L-scheme with $\mathfrak{L} = 1.5\frac{\mathcal{R}}{2} \geq \frac{\max\left( 3\mathcal{R} \cdot (p^{\varepsilon})^2\right) }{2}$.

\renewcommand{\arraystretch}{1.1}
\begin{table}[htpb!]
	\begin{center}
		\begin{tabular}{c| c  c }
			$\varepsilon$      & H        & Relative error ($e_H$) \\\hline
			$\nicefrac{1}{8}$  & $0.1768$ & $0.1938$               \\
			$\nicefrac{1}{8}$  & $0.0884$ & $0.1287$               \\
			$\nicefrac{1}{8}$  & $0.0442$ & $0.0856$               \\\hline
			$\nicefrac{1}{16}$ & $0.1768$ & $0.1797$               \\
			$\nicefrac{1}{16}$ & $0.0884$ & $0.1138$               \\
			$\nicefrac{1}{16}$ & $0.0442$ & $0.0724$               \\\hline
			$\nicefrac{1}{32}$ & $0.1768$ & $0.1690$               \\
			$\nicefrac{1}{32}$ & $0.0884$ & $0.1030$               \\
			$\nicefrac{1}{32}$ & $0.0442$ & $0.0621$               \\\hline
		\end{tabular}
	\end{center}
	\caption{History of convergence of the error for three values of $\varepsilon$ and three coarse meshes.}\label{tab:RelError_eps}
\end{table}

\Cref{tab:RelError_eps} shows the history of convergence of the error for different values of $\varepsilon$ and three coarse meshes $\trian_H$ without refinement and $H\gg h^\epsilon$.
The relative $L^2$-error $e_H$ in \Cref{tab:RelError_eps} is $e_H = \nicefrac{ \| \Pi_{h^{\varepsilon}}(p_H) - p_h \|_{L^2(\trian_{h^{\varepsilon}})}}{ \| p_h \|_{L^2(\trian_{h^{\varepsilon}})}}$ where $\Pi_{h^{\varepsilon}}(p_H)$ is the projection of the coarse-scale solution in the fine mesh $\trian_{h^{\varepsilon}}$.
With this result we show how the homogenized solution tends to the solution of the original problem when $H\to 0$ and also when $\varepsilon\to 0$.

\begin{figure}[htpb!]
	\centering
	\subfloat{\includegraphics[width=.45\textwidth]{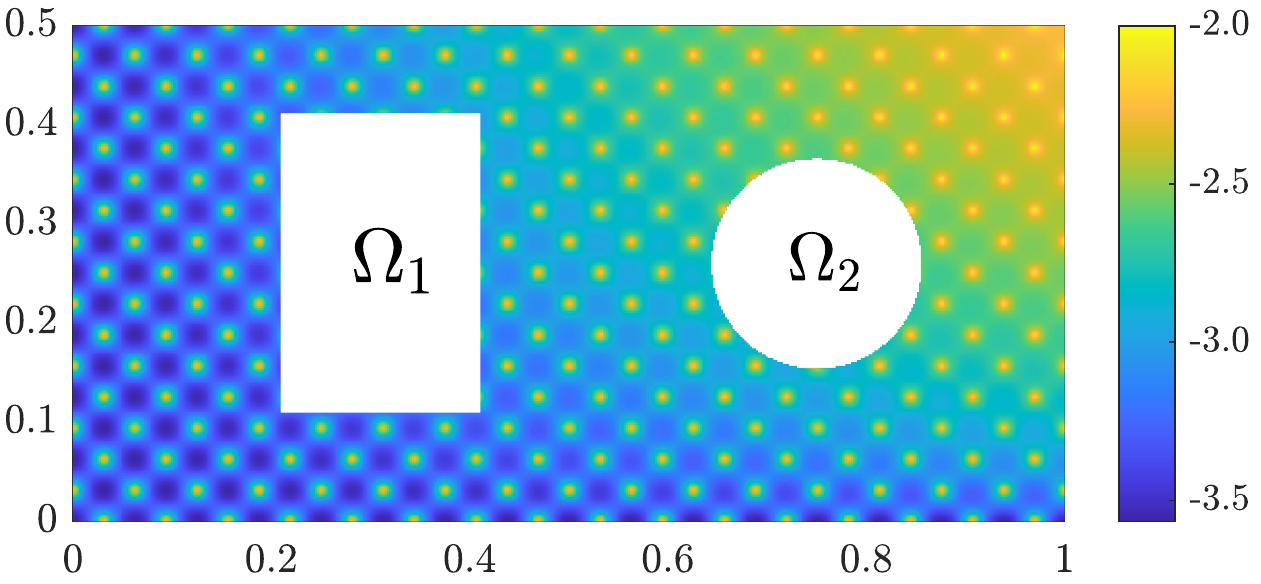}}
	\hspace{0.5cm}
	\subfloat{\includegraphics[width=.45\textwidth]{Figure5_2}}
	\caption{Fine scale permeability field ($\kep_{1,1}$) (left) $\varepsilon=\frac{1}{8}$ and (right) $\varepsilon=\frac{1}{16}$ ($Log_{10}$ scale).}\label{fig:permPeriodic}
\end{figure}

\begin{figure}[htpb!]
	\centering
	\subfloat{\includegraphics[width=0.45\textwidth]{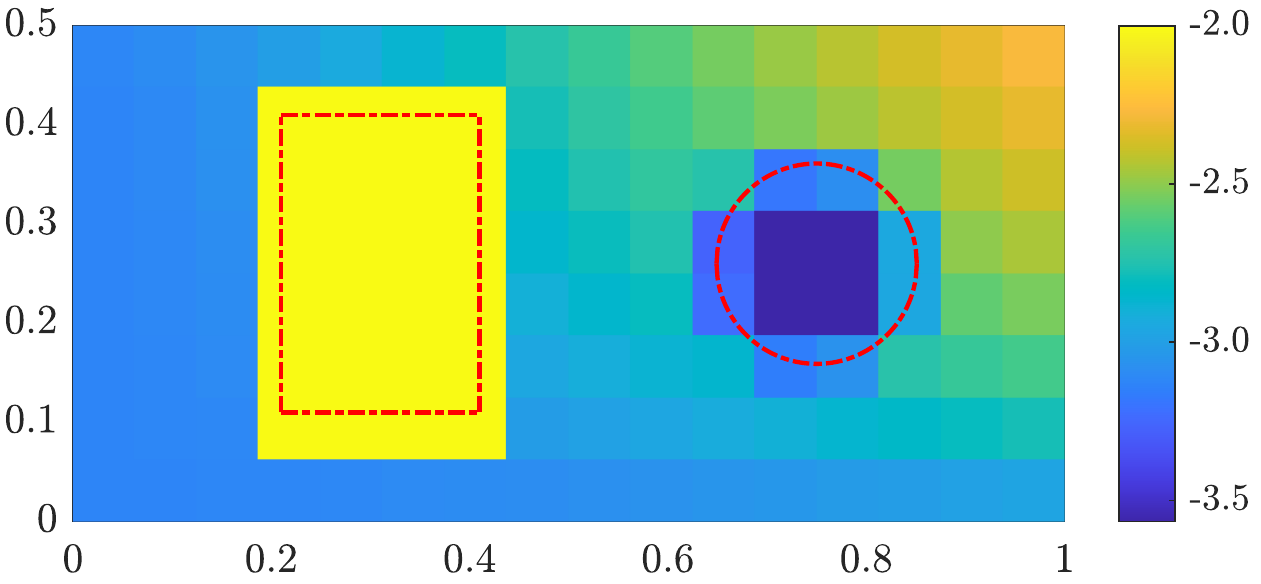}}
	\hspace{0.5cm}
	\subfloat{\includegraphics[width=0.45\textwidth]{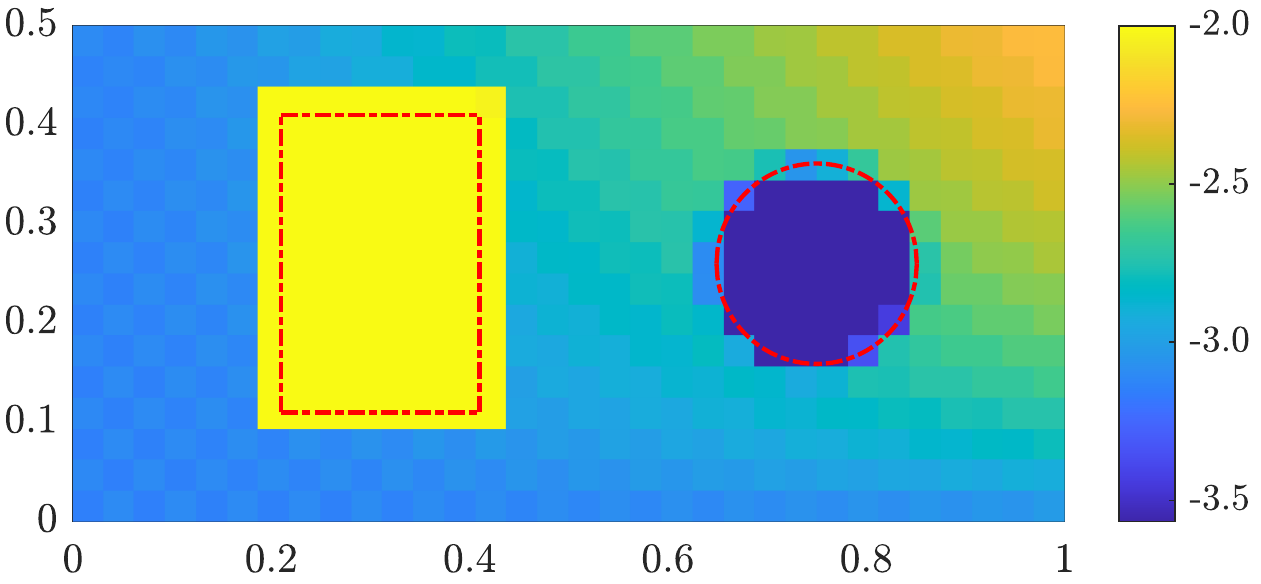}}\\
	\subfloat{\includegraphics[width=0.45\textwidth]{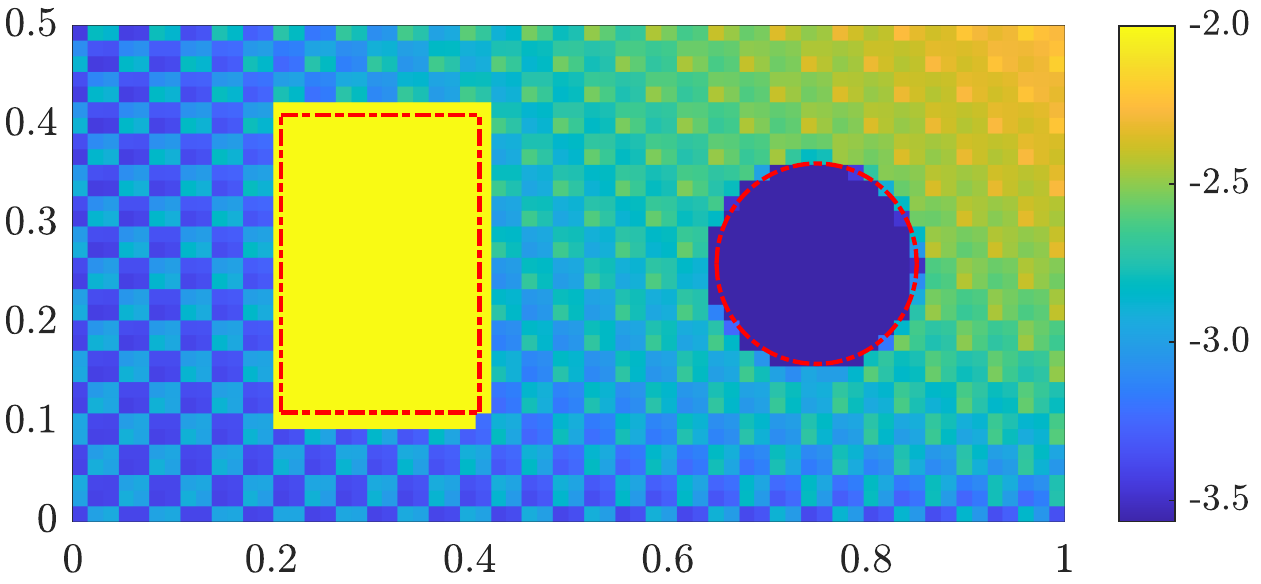}}
	\hspace{0.5cm}
	\subfloat{\includegraphics[width=0.45\textwidth]{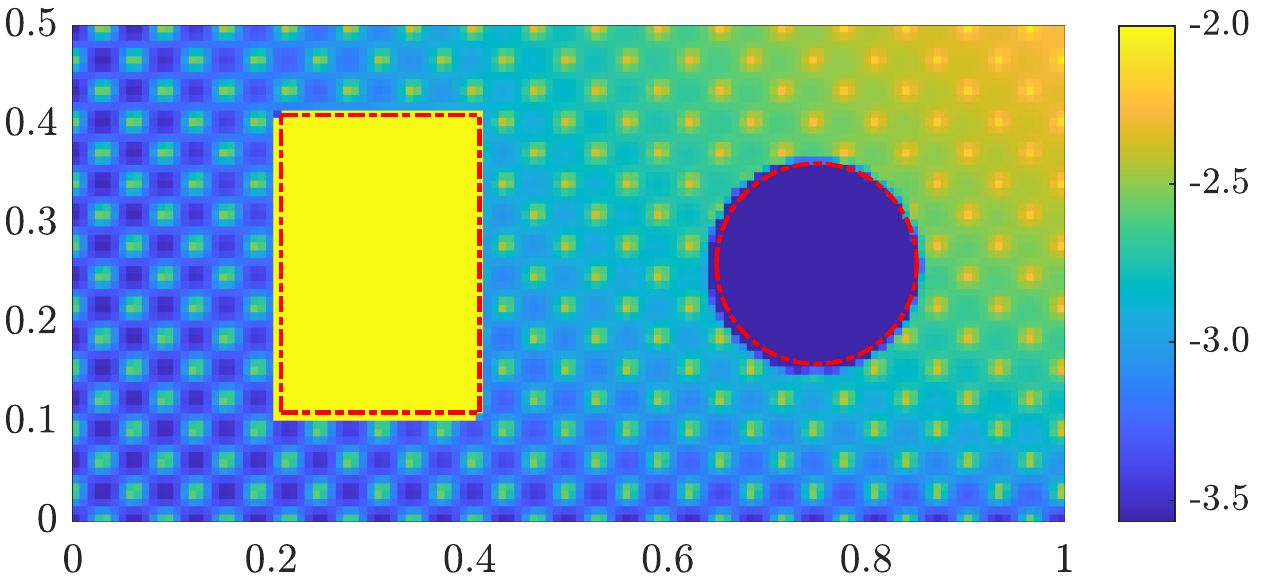}}
	\caption{Coarse-scale permeability distribution ($\kef_{1,1}$) ($Log_{10}$ scale) starting with a coarse grid of $16\times 8$ cells. The red lines indicates the original location of the low permeability zone $(\kep =10^{-7} \mathbb{I}_{2\times2})$ and high permeability zone $(\kep =10^{-2} \mathbb{I}_{2\times2})$.}\label{fig:coarseperm_per}
\end{figure}

Nevertheless, we use a modified permeability field to indicate that any assumption of periodicity is essential. We include in the same domain $\Omega$ a high permeability region $\Omega_1$ and a low permeability region $\Omega_2$ where the scalar permeability is $10^{-2}$ and $10^{-7}$ respectively.
\begin{equation*}
	\Omega_1 := [0.21,0.41]\times[0.11,0.41]  \text{ and }
	\Omega_2 := \left\lbrace \xnn\in \omep \, | \, \|\xnn-[0.75,0.26]\|_2 \leq 0.1^2 \right\rbrace.
\end{equation*}
In \Cref{fig:permPeriodic} the normalized (quasi-periodic) permeability field is shown for two values of the scale parameter $\varepsilon$. In this case the boundary conditions, the volumetric concentration, the source term and the time discretization remain the same as before.
\Cref{fig:coarseperm_per} shows four levels of the first component of the effective permeability tensor ($\kef_{1,1}$) with $\varepsilon = \frac{1}{16}$ starting with a coarse grid of $16\times 8$ cells. Referring to the different levels of the effective permeabilities, it is important to remark that the coarse-scale permeabilities are computed in zones that not always match with the initial resolution or periodicity. Here one can notice the influence of neighbouring macro-cells in the numerical solution of the micro problems $\mathbf{Ph_j}$. This effect is evident at the boundary of the low permeability zone $\Omega_2$. To point out this behaviour in the \Cref{fig:coarseperm_per} we highlight with a dashed lines the original location of the low and high permeability areas.

\begin{figure}[htpb!]
	\centering
	\includegraphics[width=.6\textwidth]{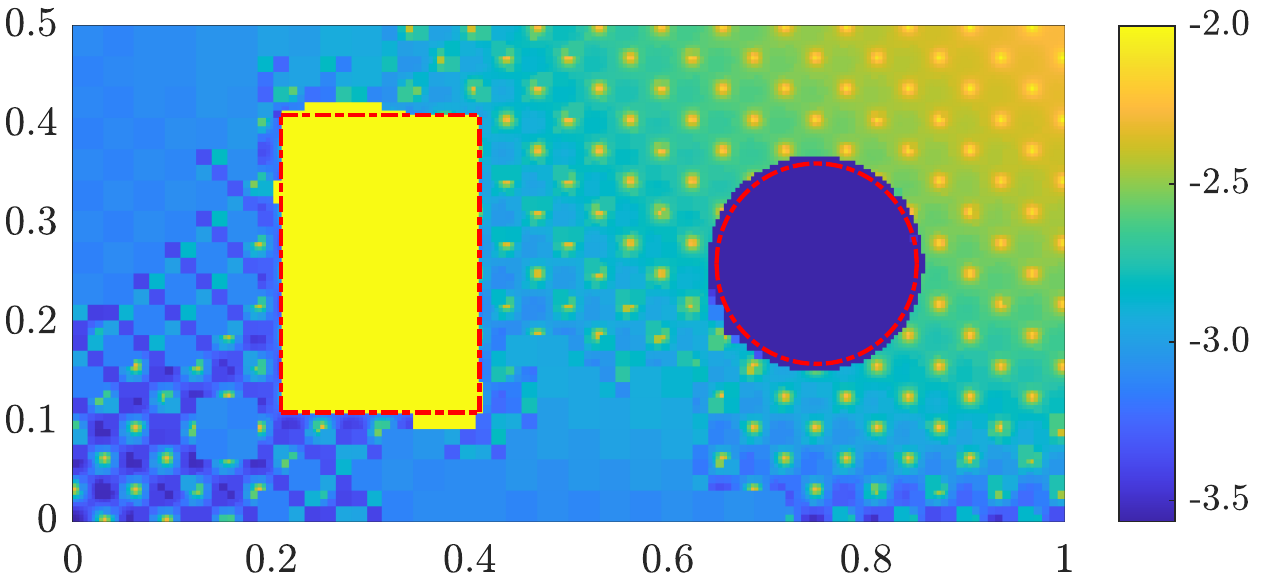}
	\caption{Refined permeability field ($\kef_{1,1}$) at $t=1$ ($Log_{10}$ scale).}\label{fig:permeabilityAdapt_perm}
\end{figure}

The anisotropic deviation of the effective permeability tensor corresponds to $9.65\cdot 10^{-5} \leq \tau_1 \leq 3.18\cdot 10^{-4}$ and $3.57\cdot10^{-5} \leq \tau_2 \leq 8.06\cdot10^{-4}$ . With this we conclude that the non-diagonal components of $\kef$ can be neglected and we use the diagonal effective tensor $\kef_D$ in  $\mathbf{PH_n^i}$. Moreover, due to the similarity between $\kef_{1,1}$ and $\kef_{2,2}$ in \Cref{fig:coarseperm_per,fig:permeabilityAdapt_perm} we only show the first component ($\kef_{1,1}$) of the effective parameter.

Furthermore, after the adaptivity process we obtain a refined version of the permeability field and \Cref{fig:permeabilityAdapt_perm} shows the result of the refined permeability at $t=1$.

\begin{figure}[htpb!] 
	\centering
	\subfloat{\includegraphics[width=.45\textwidth]{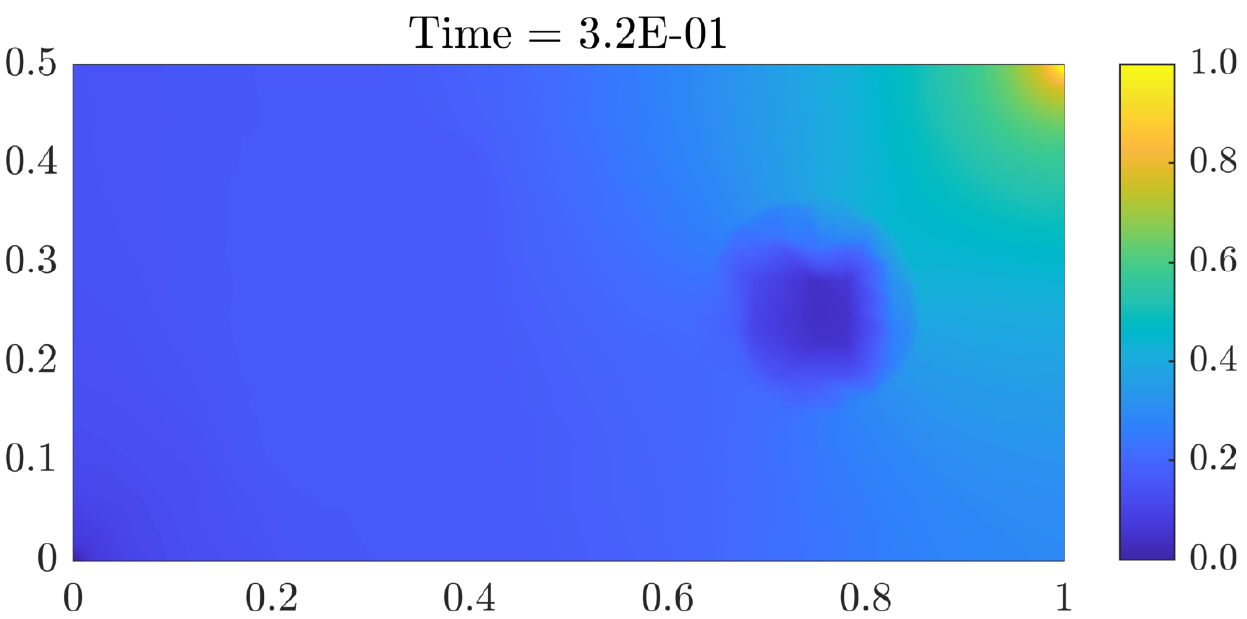}}
	\hspace{0.5cm}\subfloat{\includegraphics[width=.475\textwidth]{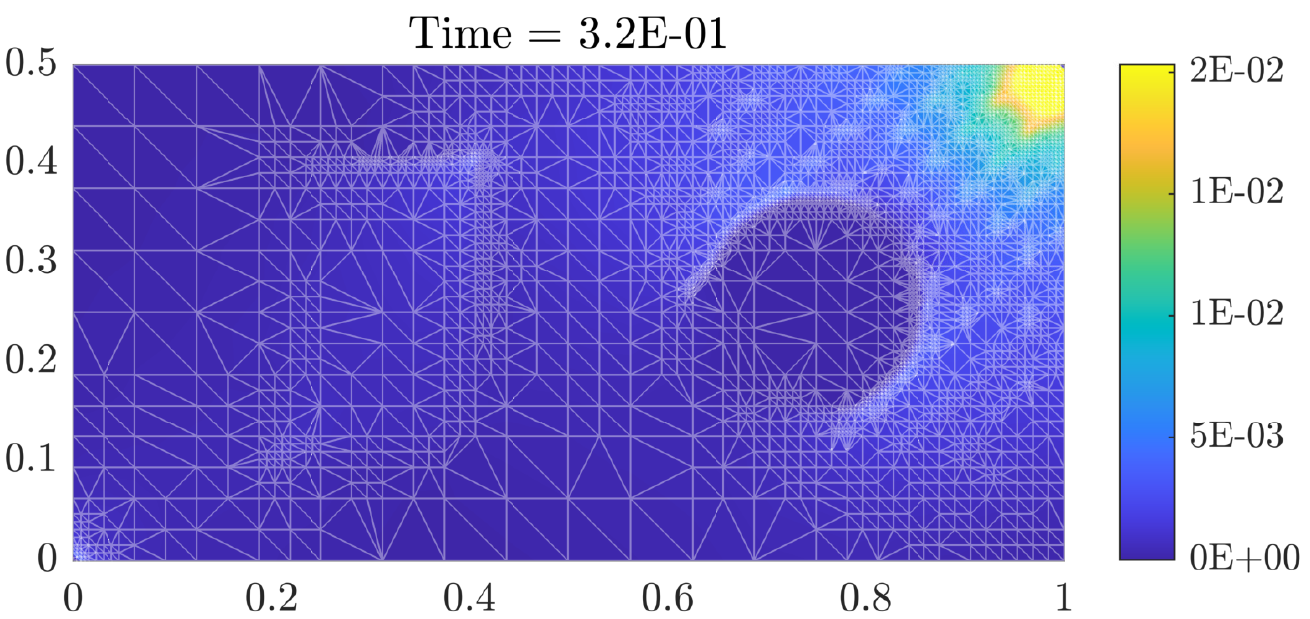}}\\
	\subfloat{\includegraphics[width=.45\textwidth]{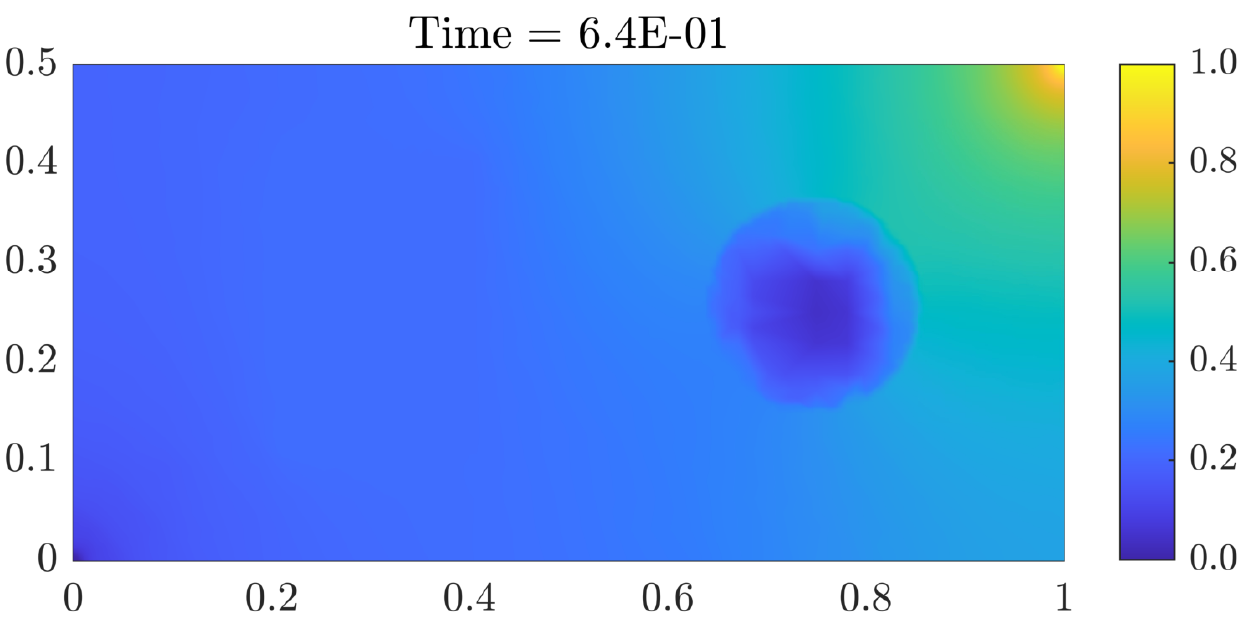}}
	\hspace{0.5cm}\subfloat{\includegraphics[width=.475\textwidth]{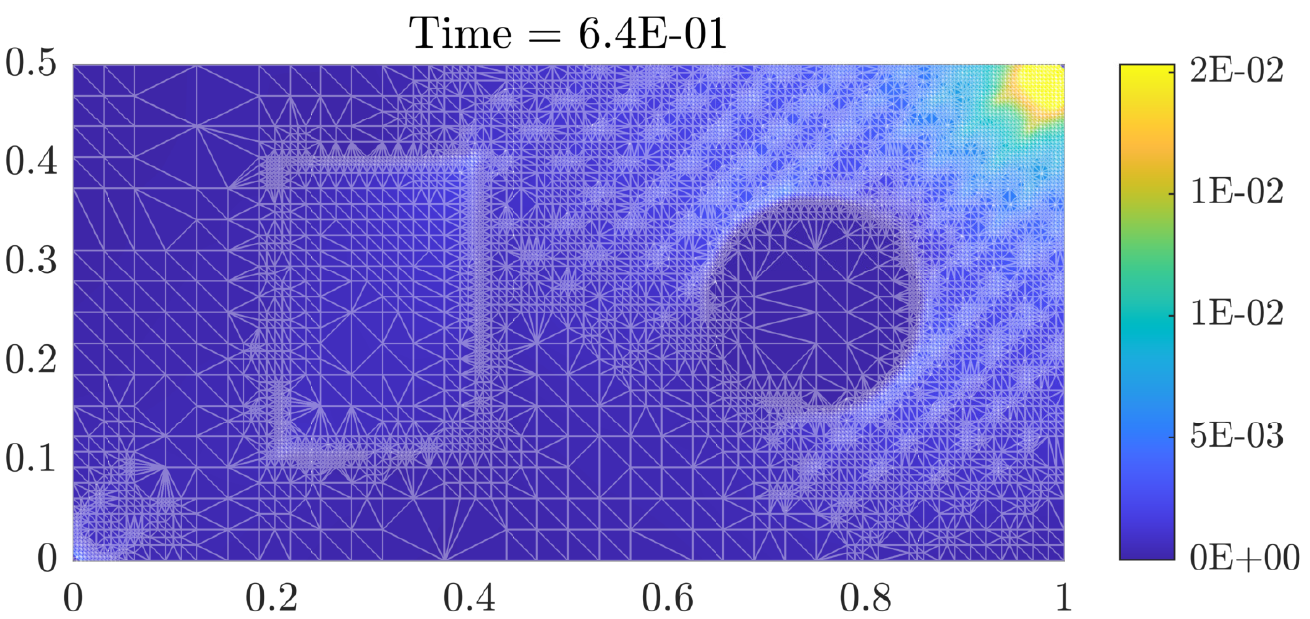}}\\
	\subfloat{\includegraphics[width=.45\textwidth]{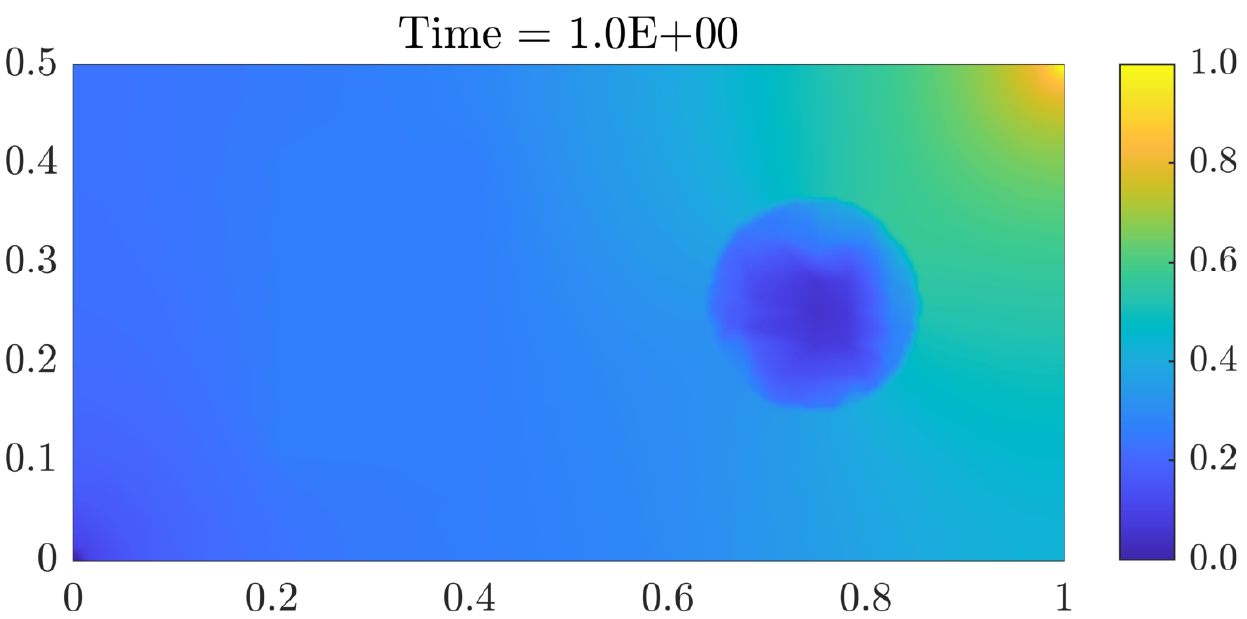}}
	\hspace{0.5cm}\subfloat{\includegraphics[width=.475\textwidth]{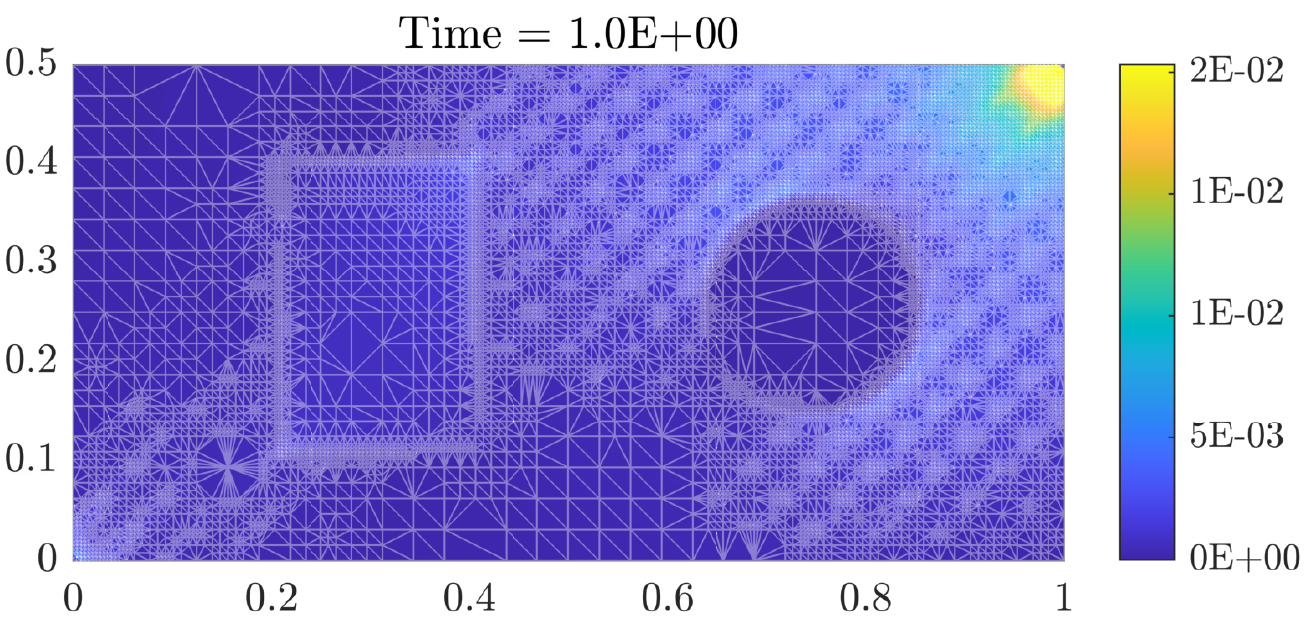}}
	\caption{Adaptive homogenization at $t=16\dt$ (top), $32\dt$ (middle), $50\dt$ (bottom). Pressure $p_{H_n}$ (left) and  magnitude of the velocity field $\|\unn_{H_n}\|_2$ (right) over meshes with $2.367$, $5.950$ and $9.659$ coarse elements. }\label{fig:SOLUTIONAdapt_per}
\end{figure}

The numerical solution of the linear upscaled problem $\mathbf{PH_n^i}$ is showed in \Cref{fig:SOLUTIONAdapt_per}. The upscaled solution is computed using the mesh refinement  described in \Cref{sec:noperiodic} using $\Theta_r = 0.5$ and $\Theta_c = 1$.
The relative $L^2$-error of the upscaled pressure $p_{H_N}$ (at the last time step) is calculated as follows
\begin{equation}
	E_\mathrm{T}^2 = \frac{\int_0^\mathrm{T} \|p_{H_N} - p_{h^\varepsilon}\|_\Omega^2 dt}{\int_0^\mathrm{T} \|p_{h^\varepsilon}\|_\Omega^2 dt} = 0.0507
\end{equation}
Here, the $L^2$-error of the upscaled pressure $p_{H_N}$ is $E_\mathrm{T}^2 = 1.6\%$ using only $14.7\%$ of the degrees of freedom used to compute a fine scale solution with mesh size $h^\varepsilon\ll\varepsilon = \frac{1}{16}$.

\begin{figure}[htpb!]
	\centering
	\includegraphics[width=0.5\textwidth]{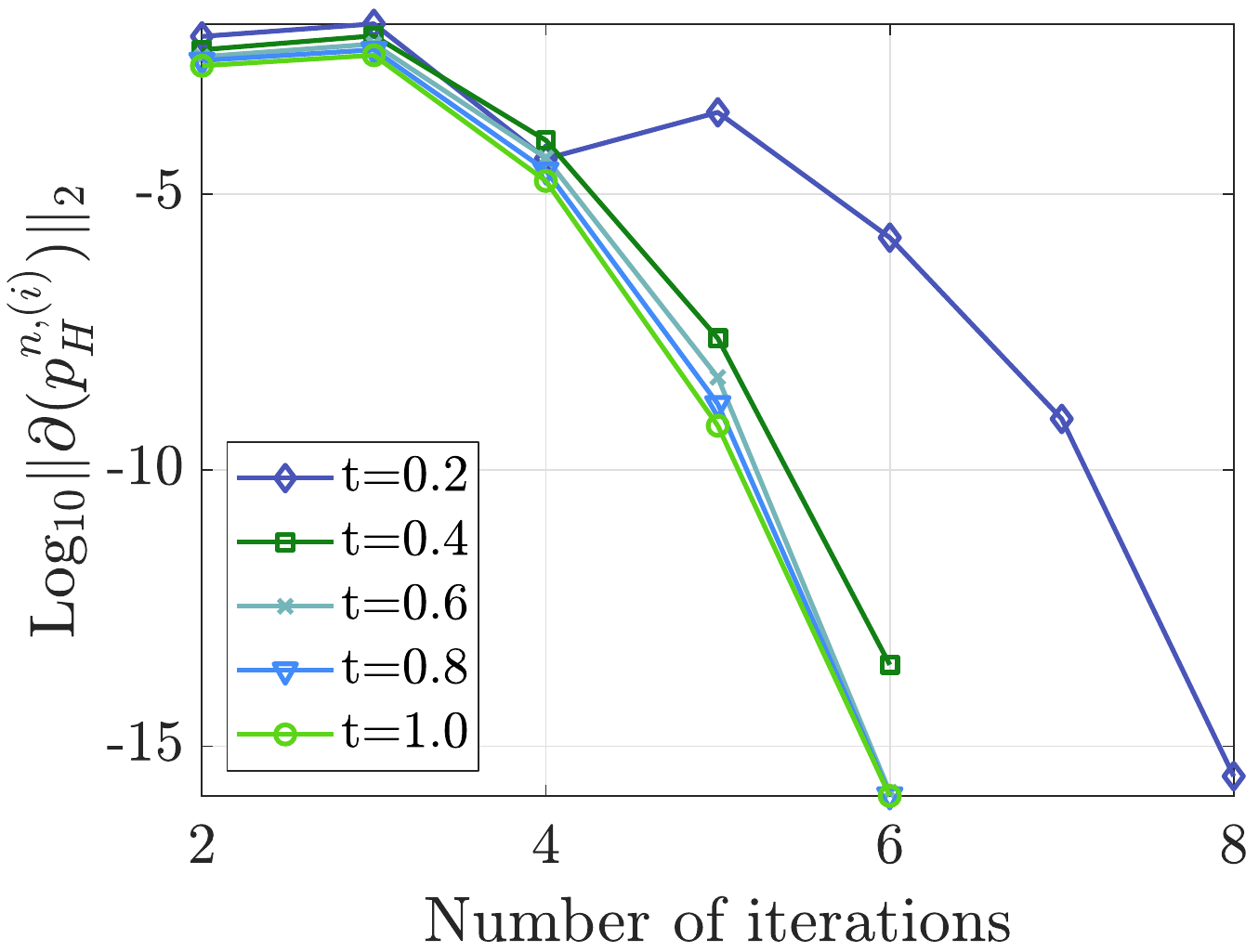}
	\caption{Convergence of the residual in the non-linear solver. Results for four different times steps using the L-scheme with $\mathfrak{L}=1.5\frac{\mathcal{R}}{2}$ and the Newton method afterwards.}\label{fig:convergece_Adapt2_perm}
\end{figure}

Concerning the behaviour of the non-linear solver, our test case is an example where the convergence of the Newton method highly depends on the initial guess. However the convergence of the L-scheme is not optimal; i.e., even though the L-scheme converges we do not want to lose the quadratic convergence of the Newton method.
To compute the solution in \Cref{fig:SOLUTIONAdapt_per}, the L-scheme reaches the threshold $\|\partial(p_H^{n,(i)})\|_2 < 10^{-10}$ after an average of $70$ iterations. In order to improve the linear solver we use a mixed strategy (see \cite{list2016study}). The target is to construct an initial solution that suits a non-problematic starting point for the Newton method. In this case we used the L-scheme until $\|\partial(p_H^{n,(i)})\|_2 < 10^{-2}$ and then the classical Newton method until one reaches  $\|\partial(p_H^{n,(i)})\|_2 < 10^{-10}$ (see \Cref{fig:convergece_Adapt2_perm}).

\subsection{The non-periodic case}
\label{sec:Num2}
\begin{figure}[htpb]
	\centering
	\includegraphics[width=0.65\textwidth]{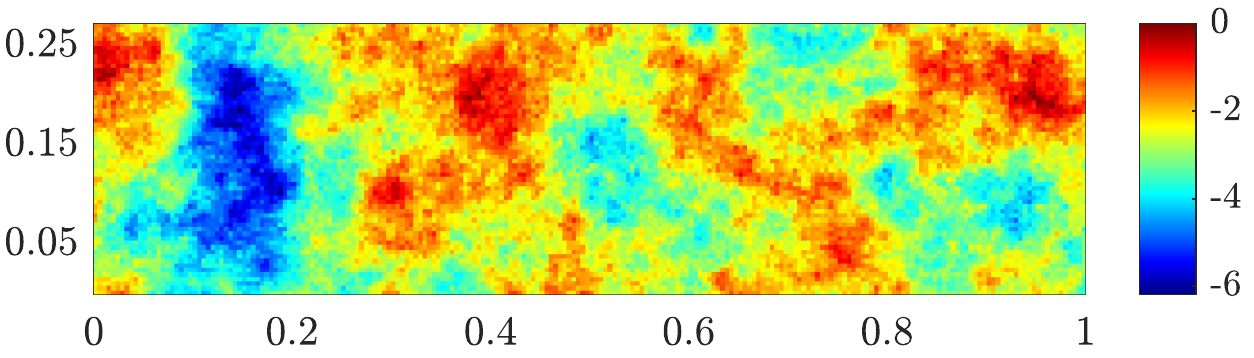}
	\caption{Fine scale permeability distribution ($\kep_{1,1}$) for SPE10th-TopLayer  ($Log_{10}$ scale).}\label{fig:permspe10}
\end{figure}

Here we consider a highly heterogeneous and non-periodic medium. We utilize the data of the SPE Comparative Solution Projects \cite{christie_tenth_2001}. This provides a vehicle for independent comparison of methods and a recognized suite of test datasets for specific problems. Our isotropic permeability field $\kep$ is defined by the top field of SPE10th data set (see \Cref{fig:permspe10}).

The macro-scale domain is a two-dimensional rectangle (see \Cref{fig:permspe10}). External boundaries are impermeable; i.e., we take no-flux boundary conditions. The domain is initialized with pressure $p_0 = 0$. A source and a sink are placed in the lower-left and the upper-right corners, having fixed pressures of $1$ and $0$, respectively.

Moreover, the volumetric concentration is $b^{\varepsilon}(\xnn,(p^{\varepsilon})) = \mathcal{R}\cdot(p^{\varepsilon})^3$. Here $\mathcal{R}$ is defined as in \cref{sec:Num1}. For the time discretization we take $\mathrm{T}=1$ with $\dt=0.02$ and the parameter for the non-lineal solver is $\mathfrak{L} = 1.5\frac{\mathcal{R}}{2} \geq \frac{\max\left( 3\mathcal{R}\cdot(p^{\varepsilon})^2\right) }{2}$. The criteria for the dynamic mesh refinement, described in \Cref{sec:noperiodic}, are $\Theta_r=0.2$ and $\Theta_c=10$.

\begin{figure}[htpb!] 
	\centering
	\subfloat{\includegraphics[width=.48\textwidth]{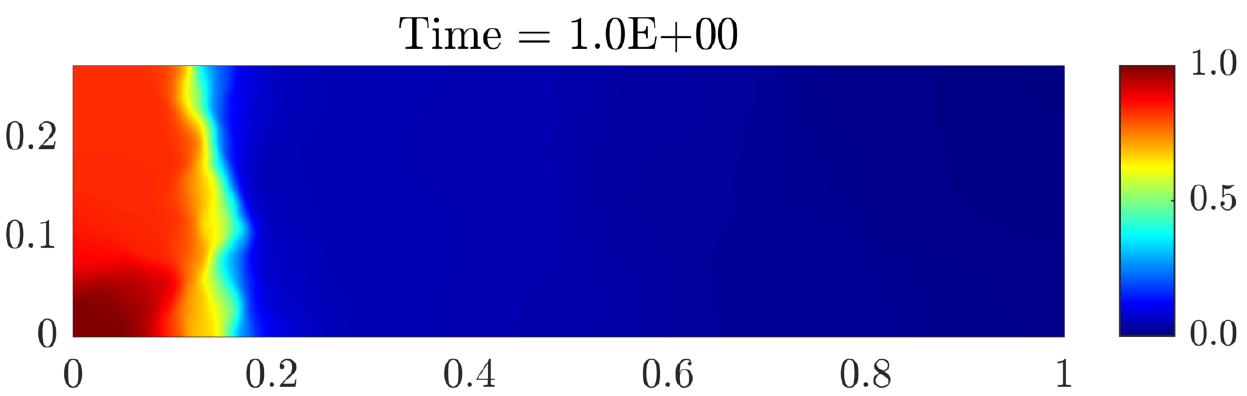}}
	\hspace{0.2cm}
	\subfloat{\includegraphics[width=.49\textwidth]{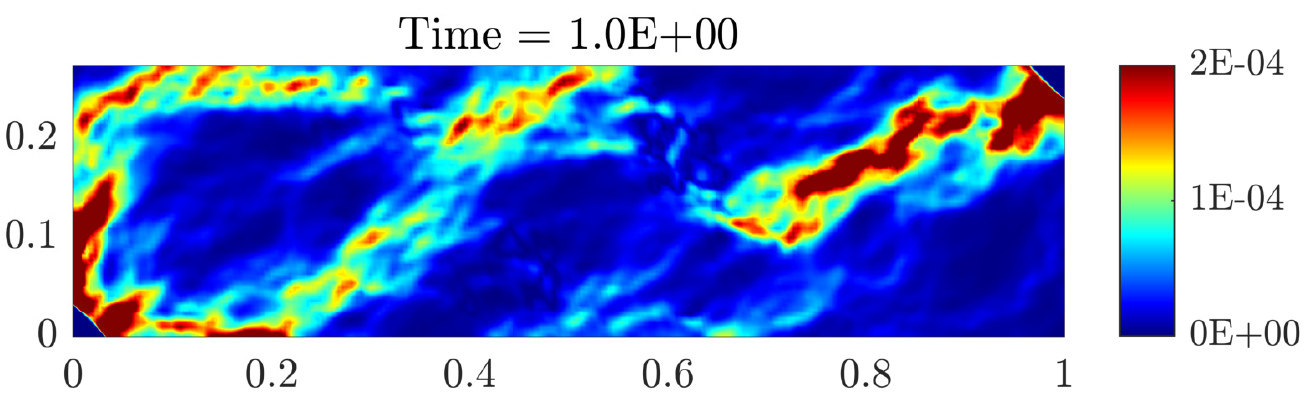}}
	\caption{Fine scale pressure $p_h^\varepsilon$ (left) and (right) magnitude of the velocity field $\|\unn_h^\varepsilon\|_2$.}
	\label{fig:SOLUTIONref}
\end{figure}

\begin{figure}[htpb!]
	\centering
	\includegraphics[width=0.65\textwidth]{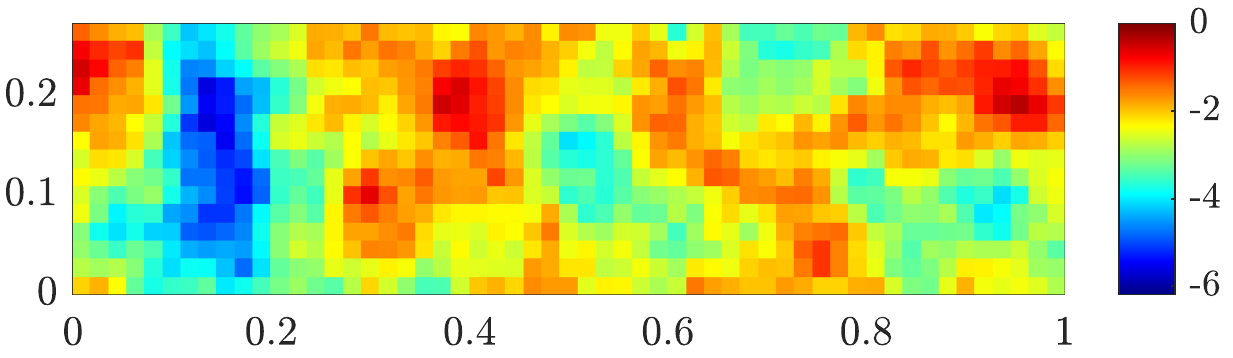}
	\caption{Coarse-scale permeability distribution ($\kef_{1,1}$) ($Log_{10}$ scale).}\label{fig:coarseperm}
\end{figure}

\begin{figure}[htpb!]
	\centering
	\includegraphics[width=0.65\textwidth]{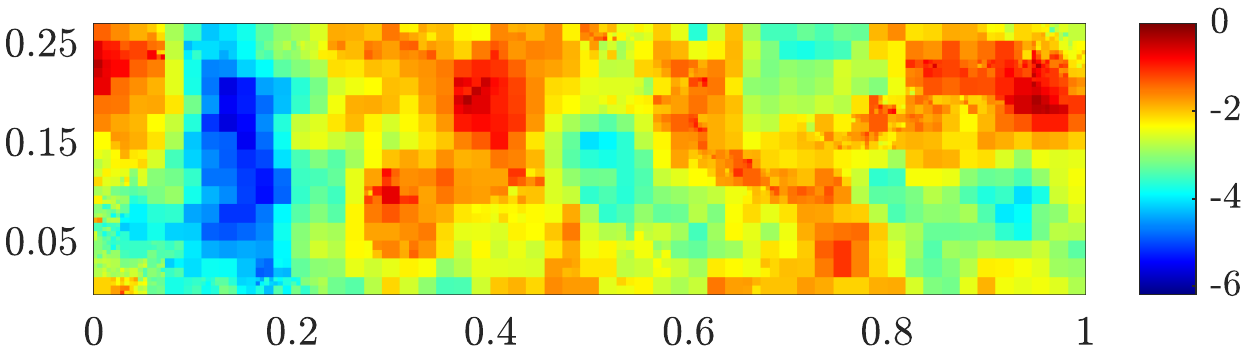}
	\caption{Refined permeability field ($\kef_{1,1}$) at $t=1$  ($Log_{10}$ scale).}\label{fig:permeabilityAdapt}
\end{figure}

To solve the problem \eqref{eq:mixedOrig} with the resolution of \Cref{fig:permspe10} we construct a grid with $26.400$ elements in a homogeneous triangular mesh $\trian_{h^\varepsilon}$. In \Cref{fig:SOLUTIONref} we show the reference solution $(p_{h^\varepsilon},\unn_{h^\varepsilon})$ at the last time step.
Using a coarse grid of $55\times 15$ squares we compute the first effective permeability field. This coarse grid corresponds to a macro-scale mesh with $1.650$ triangular elements. In \Cref{fig:coarseperm} we show the first component ($\kef_{1,1}$) of the coarse-scale permeability field. The anisotropic deviation encountered in this effective permeability $\kef$ corresponds to $\tau_1 = 2.5\cdot 10^{-4}$ and $\tau_2 = 1.7\cdot 10^{-3}$. For this reason, the non-diagonal components of $\kef$ are neglected and in \Cref{fig:coarseperm,fig:permeabilityAdapt} we display only the first component ($\kef_{1,1}$) of the effective tensor .

\begin{figure}[htpb!]
	\centering
	\includegraphics[width=0.65\textwidth]{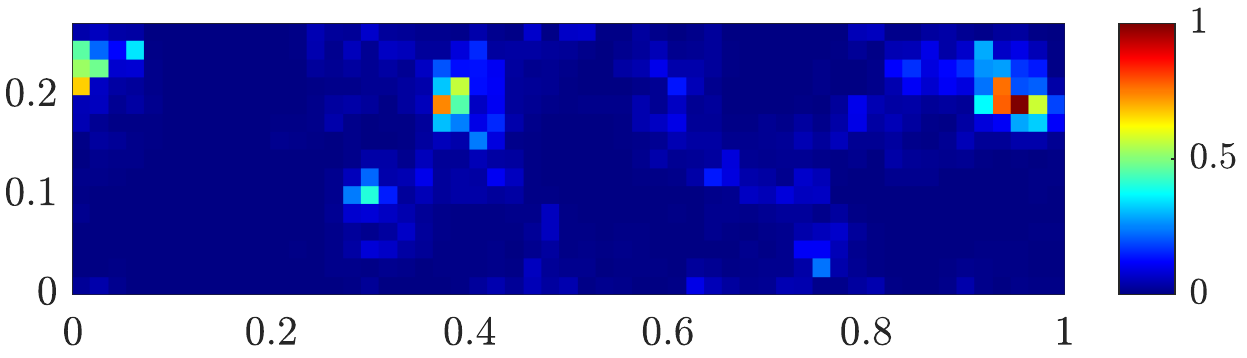}
	\caption{Difference between the coarse-scale effective permeabilities using homogenization vs harmonic average.}\label{fig:coarsepermdiff}
\end{figure}

In \Cref{fig:coarsepermdiff} we show the difference between the effective permeabilities computed with homogenization and using the harmonic average. The difference between these strategies is higher in zones with high permeability and one can point out that the harmonic average always underestimates the permeability. This is problematic because the high permeability regions are regions where one should increase the accuracy of the effective parameter in order to have better numerical solutions. When we compute the numerical solution of the problem $\mathbf{PH_n^i}$ using the harmonic average of the permeability the relative $L^2-$error of the pressure is $7\%$ higher than the error when using the effective permeability \eqref{eq:EffectiveTensor}.

Using the adaptivity process we obtain a refined version of the permeability field. \Cref{fig:permeabilityAdapt} shows the result of the permeability field after the mesh adaptivity process.

\begin{figure}[htpb!] 
	\centering
	\subfloat{\includegraphics[width=.48\textwidth]{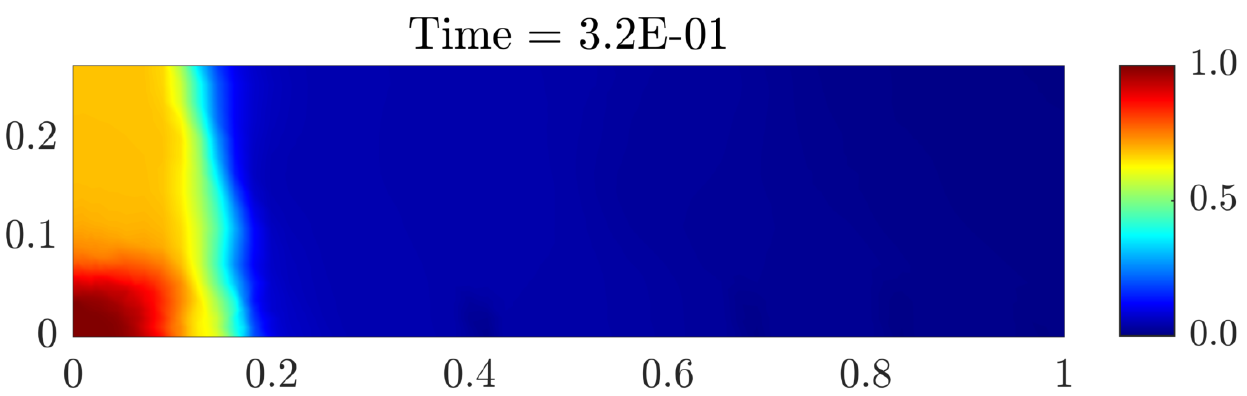}}
	\hspace{0.2cm}
	\subfloat{\includegraphics[width=.48\textwidth]{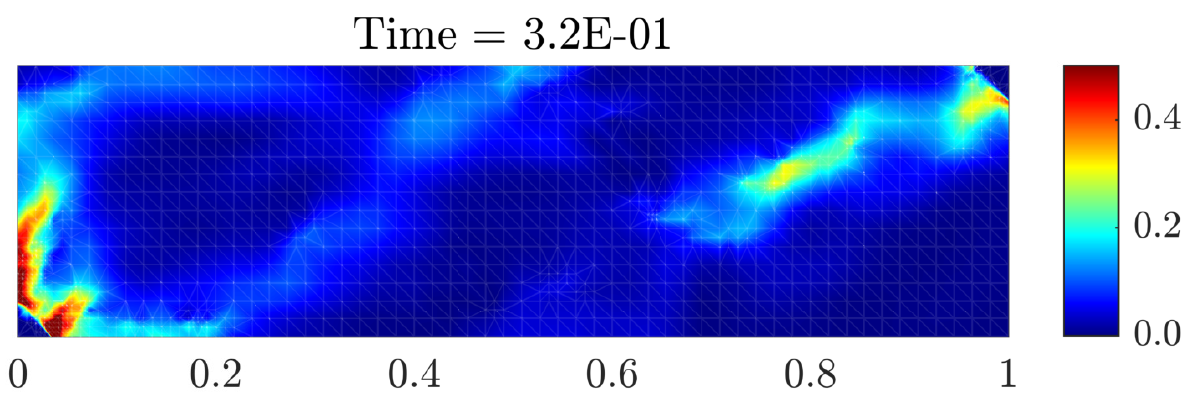}}\\
	\subfloat{\includegraphics[width=.48\textwidth]{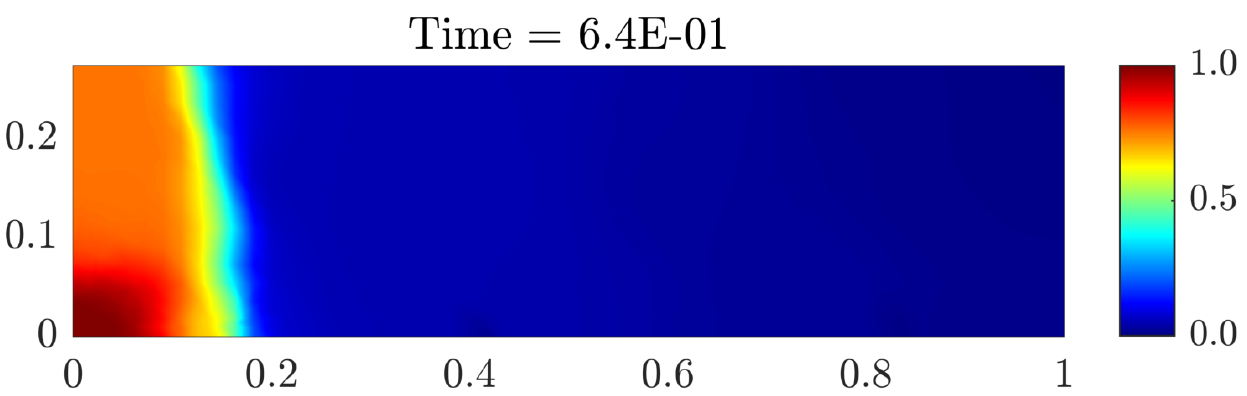}}
	\hspace{0.2cm}
	\subfloat{\includegraphics[width=.48\textwidth]{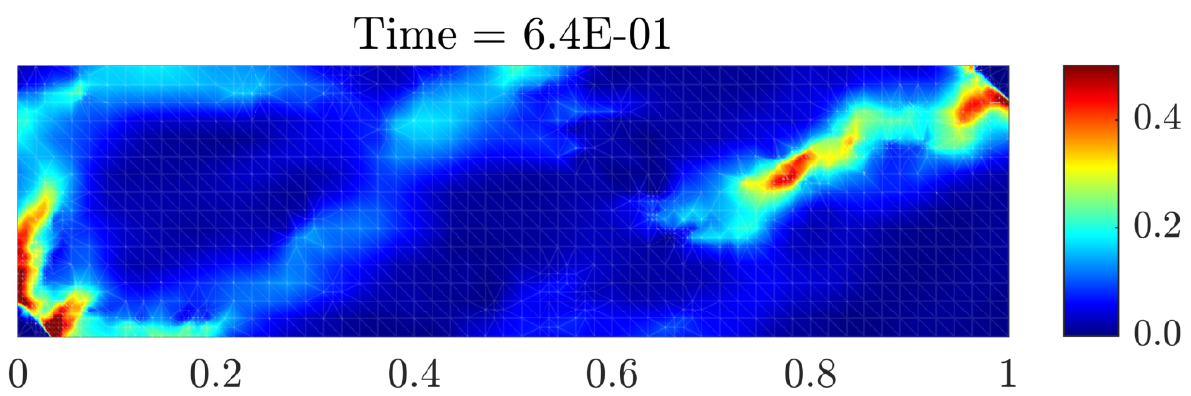}}\\
	\subfloat{\includegraphics[width=.48\textwidth]{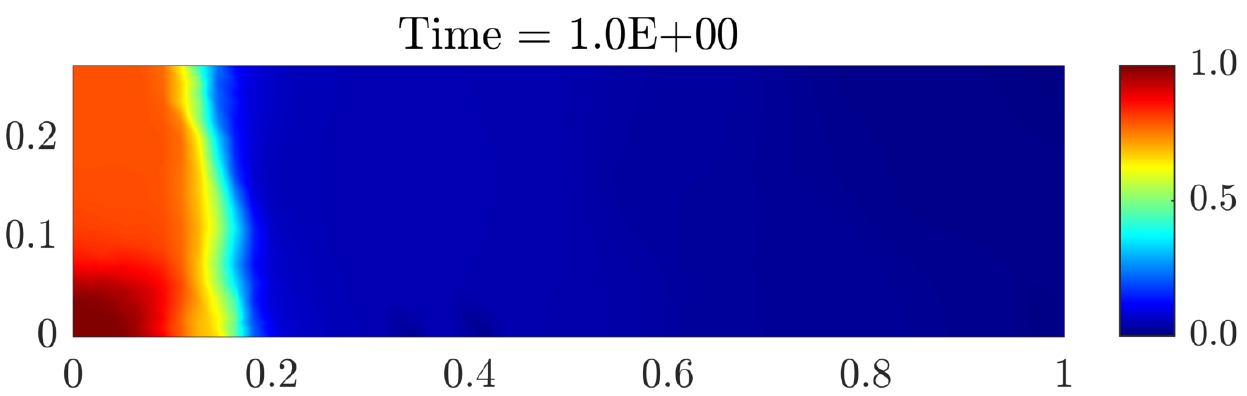}}
	\hspace{0.2cm}
	\subfloat{\includegraphics[width=.48\textwidth]{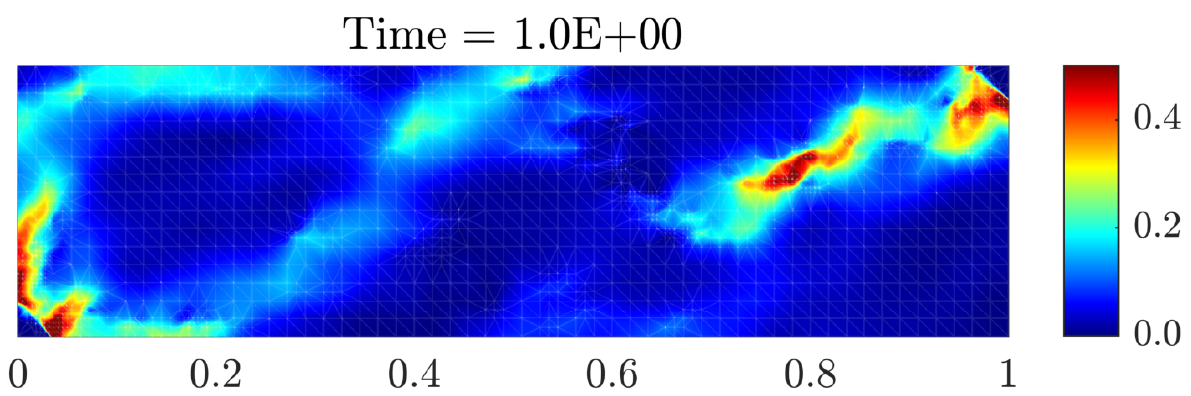}}
	\caption{Adaptive homogenization at $t=16\dt$ (top), $32\dt$ (middle), $50\dt$ (bottom). Pressure $p_{H_n}$ (left) and (right) magnitude of the velocity field $\|\unn_{H_n}\|_2$ over meshes with $2.710$, $3.152$ and $3.783$ coarse elements. }\label{fig:SOLUTIONAdapt}
\end{figure}

\Cref{fig:SOLUTIONAdapt} shows the numerical solution of the upscaled problem $\mathbf{PH_n^i}$ using the mesh adaptive re-meshing described in \cref{sec:noperiodic}. At the end of the adaptivity process, the relative $L^2$-error of the upscaled pressure $p_{H}$ is $E_\mathrm{T}^2 = 5.07\%$ using only the $14.5\%$ of degrees of freedom used to computed the reference solution $p_{h^\varepsilon}$.

\begin{figure}[htpb!]
	\centering
	\includegraphics[width=0.5\textwidth]{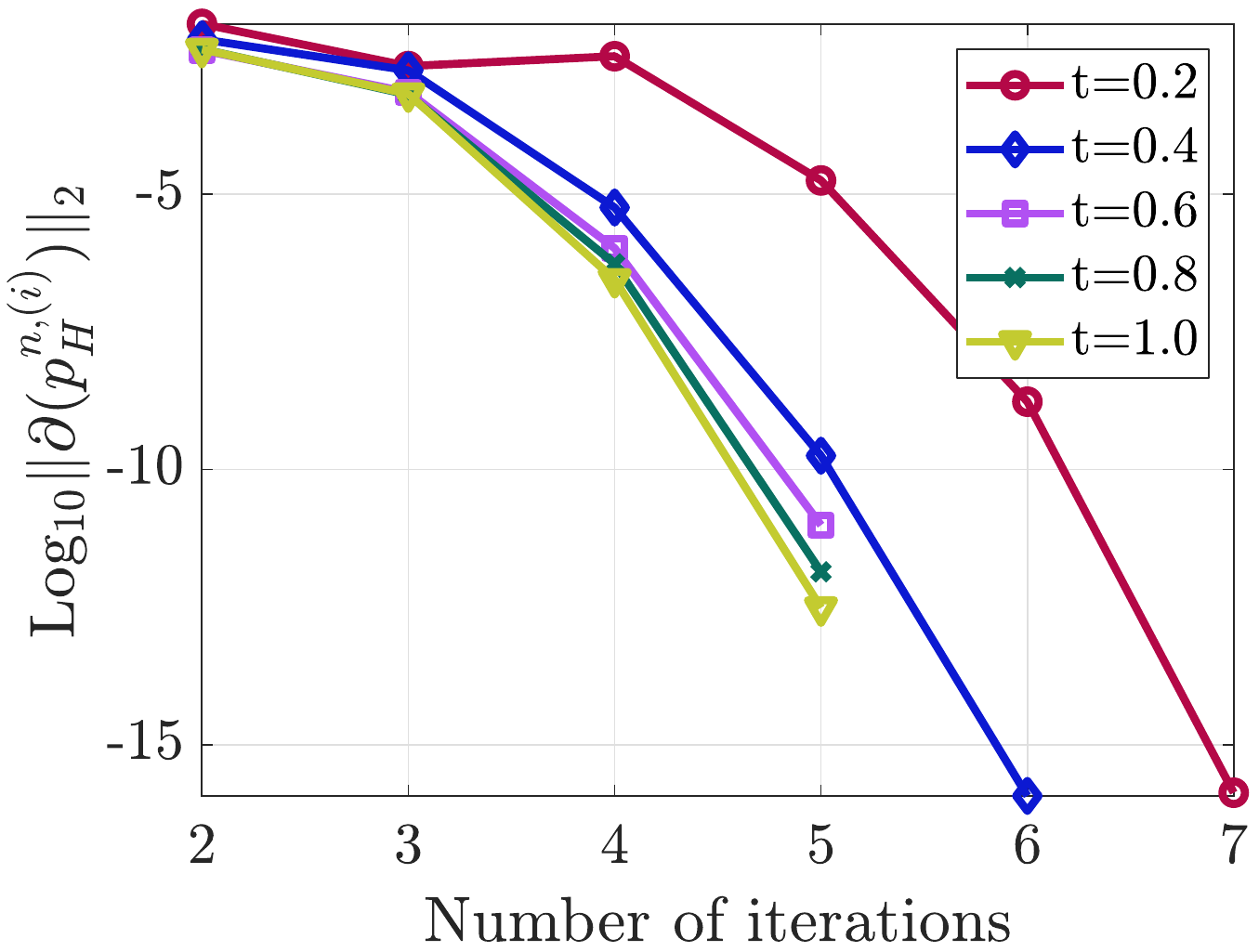}
	\caption{Convergence of the residual in the non-linear solver. Results for five different times steps using the L-scheme with $\mathfrak{L}=1.5\frac{\mathcal{R}}{2}$ and the Newton method afterwards.}\label{fig:convergece_Adapt2}
\end{figure}

Finally, in \Cref{fig:convergece_Adapt2} we show the convergence of the norm of the residual $\partial(p_H^{n,(i)})$ when one uses a combination of the L-scheme and Newton method. Here we use a mixed strategy (see \cite{list2016study}) to construct an initial solution that suits a non-problematic starting point for the Newton method. In this case we use the L-scheme until $\|\partial(p_H^{n,(i)})\|_2 < 10^{-2}$ and then the classical Newton method until one reaches  $\|\partial(p_H^{n,(i)})\|_2 < 10^{-10}$ and as we see in \Cref{fig:convergece_Adapt2} the quadratic convergence of the Newton method is recovered.

\section{Conclusions}
\label{sec:conclusions}

We have presented a numerical scheme based on homogenization to solve a non-linear parabolic equation with highly oscillatory characteristics. The discrete non-linear system is obtained by a  backward Euler and lowest order Raviart-Thomas mixed finite element discretization.
Our approach proposes a local mesh adaptivity that leads to the computation of the effective parameters locally through decoupled cell problems. The mesh adaptivity is based on the idea that the upscaled parameters are updated only when it is necessary. Moreover, to illustrate the performance we have presented two general examples. We constructed a periodic case to show the history of convergence of the error when the scale separation tends to zero. In the non-periodic case we used a benchmark from the SPE10th project and we showed that the homogenization can be used in more general non-periodic cases.

In addition, we combined the standard Newton method and the L-scheme to improve the behaviour of the non-linear solvers. We presented a combination of techniques that led to a very efficient numerical scheme. It is relevant to mention that besides the theory mentioned in this paper the applicability of this strategy is vast. Extensions of our adaptive algorithm including more complex micro-scale models are applicable.  Those include from reactive transport up to moving interfaces affecting the structure of the micro-scale.

\section*{Acknowledgements}
The authors gratefully acknowledge financial support from the Research Foundation - Flanders (FWO) through the Odysseus programme (Project G0G1316N). In addition, we wish to thank Professor Mary F. Wheeler and Professor Ivan Yotov who made valuable suggestions or who have otherwise contributed to the ideas behind this manuscript. Part of this work was elaborated during the stay of the first author in the University of Bergen supported by the Research Foundation - Flanders (FWO) through a travel grant for a short stay abroad.

\bibliographystyle{elsarticle-num}
\bibliography{biblio0}
\end{document}